\documentclass{article}
\usepackage{epsfig}
\title{Locally constant $n$-operads as higher braided operads}
\author{M.A. Batanin\protect \footnote{The author holds the Scott Russell Johnson Fellowship in
the Centre of Australian Category Theory at Macquarie University}\\ Macquarie University,  NSW 2109, Australia \\
e-mail: {mbatanin@ics.mq.edu.au}}

\newtheorem{theorem}{\bf Theorem}[section]
\newtheorem{defin}{\bf Definition}[section]
\newtheorem{pro}{\bf Proposition}[section]
\newtheorem{lem}{\bf Lemma}[section]

\newcommand{\Sm}{\mbox{${\cal S}$}}
\newcommand{\Br}{\mbox{${\cal B}r$}}

\newcommand{\W}{\mbox{$\cal W$}}
\newcommand{\Conf}{\mbox{ Conf}}
\newcommand{\Proof}{\noindent {\bf Proof.}\hspace{2mm}} 
\newcommand{\Remark}{\noindent {\bf Remark.}\hspace{2mm}} 
\newcommand{\ph}{\mbox{$\bf ph$}}
\newcommand{\rh}{\mbox{$\bf rh$}}

\newcommand{\Q}{
{\unitlength=0.25mm
\begin{picture}(500,10)(-10,0)
\put(440,10){\line(0,-1){10}}
\put(440,0){\line(1,0){10}}
\put(450,0){\line(0,1){10}}
\put(450,10){\line(-1,0){10}}
\put(451,11){\line(-1,0){10}}
\put(451,1){\line(0,1){10}}
\put(450.5,0.5){\line(0,1){10}}
\put(450.5,10.5){\line(-1,0){10}}
\end{picture}}}

\begin{document}

\maketitle
\begin{abstract}

We introduce a category of {\it locally constant $n$-operads}   which can be considered  as  the category of  higher braided operads.   For $n=1,2,\infty$ the homotopy  category of locally constant $n$-operads is equivalent to the homotopy category of classical        
nonsymmetric, braided and symmetric operads correspondingly.

\

1991 Math. Subj. Class.  18D20 , 18D50, 55P48     
\end{abstract}

\section{Introduction}

It is well known  that  contractible nonsymmetric operads detect $1$-fold loop spaces,  contractible  
braided operads detect $2$-fold loop spaces and that  contractible  
symmetric operads detect $\infty$-fold loop spaces. A natural question  arises : is there a sequence of groups $G^{(n)}=\{G^{(n)}_k\}_{k\ge 0}$ together with a notion of $G^{(n)}$-operad, which we would call {\it $n$-braided operad}, such that the algebras of a contractible such operad are $n$-fold loop spaces? With some natural  minor assumptions  one can prove that the answer to the above question is negative. This is because for such an operad $A$ the quotient $A_k/G^{(n)}_k$ is a $K(G^{(n)}_k,1)$-space. One can show, however, that such a quotient must have a homotopy type of the space of unordered configurations of $k$ points in $\Re^n ,$ which is a $K(\pi,1)$-space only for $n=1,2,\infty .$ 

In this paper we show that there is a category of operads which we can think of as a correct replacement for the nonexistent category of  $G^{(n)}$-operads in all dimensions. We call them {\it locally constant $n$-operads.} For $n=1,2,\infty$ the homotopy  category of locally constant $n$-operads is equivalent to the homotopy category of classical        
nonsymmetric, braided and symmetric operads correspondingly.

Here is a brief overview of the paper. In section \ref{symgroup} we recall the definitions of symmetric and braided operads. In Section \ref{nord} we introduce the category of $n$-ordinals as higher dimensional analogue of  
 the category of finite ordinals. Using this category and its subcategory of quasibijections we 
define $n$-operads and quasisymmetric $n$-operads in Section \ref{qoperad}. In Section \ref{qconf} we show  that  the   category of quasibijections   is closely related to the classical Fox-Newirth stratification of configuration spaces. As a corollary we observe that the nerve of this category has homotopy type of unordered configurations of points in $\Re^n.$  We also prove two technical lemmas which we use in  Section \ref{vs} to relate different operadic notions. Finally in Section \ref{lcoperad} we introduce  locally constant operads and compare them  with symmetric, braided and quasisymmetric operads. We also state our recognition principle for $n$-fold loop spaces.

\section{Symmetric and braided operads}\label{symgroup}

 For a natural number $n$ we will denote by $[n]$ the ordinal
$$0 \ < \ 1 \ < \ \ldots \ < n .$$
We  denote an empty ordinal by $[-1].$  A morphism from $[n]\rightarrow [k]$ 
is any function between underlying sets. It can be order preserving or not. 
It is clear that we then have a category. We denote this category by $\Omega^s$. Of course,
$\Omega^s$ is equivalent to the category of finite sets.    In particular, the symmetric group
$S_{n+1}$ is the group of automorphisms of $[n]$.  

Let $\sigma:[n]\rightarrow [k]$ be a morphism in $\Omega^s$ and 
let  $ 0\le i \le k .$ Then the preimage $\sigma^{-1}(i)$ has a   linear order induced from $[n] .$ Hence, there exists a unique object   $[n_i]\in \Omega^s$ and a unique order preserving bijection $[n_i]\rightarrow \sigma^{-1}(i) .$ We will call $[n_i]$ 
the {\it fiber} of $\sigma$ over $i$ and will denote it $\sigma^{-1}(i) $ or $[n_i].$   

Analogously, given a composite of morphisms  in $\Omega^s :$ 
\begin{equation}\label{nlk} [n]\stackrel{\sigma}{\longrightarrow} [l] \stackrel{\omega}{\longrightarrow} [k] \end{equation}
we will denote
 $\sigma_i$   the {\it $i$-th fiber} of $\sigma$; i.e. the
pullback 

{\unitlength=0.9mm

\begin{picture}(60,40)(-18,-3)

\put(10,25){\makebox(0,0){\mbox{$ \sigma^{-1}(\omega^{-1}(i))$}}}
\put(10,20){\vector(0,-1){10}}
\put(12,15){\shortstack{\mbox{$ $}}}

\put(22,25){\vector(1,0){10}}

\put(24,26){\shortstack{\mbox{$\sigma_i $}}}

\put(45,25){\makebox(0,0){\mbox{$ \omega^{-1}(i)$}}}
\put(45,20){\vector(0,-1){10}}

\put(75,25){\makebox(0,0){\mbox{$[1] $}}}
\put(75,20){\vector(0,-1){10}}

\put(57,25){\vector(1,0){10}}

\put(57,21){\shortstack{\mbox{$ $}}}

\put(57,28){\shortstack{\mbox{\small $ $}}}

\put(77,14){\shortstack{\mbox{\small $\xi_i $}}}

\put(10,5){\makebox(0,0){\mbox{$ [n]$}}}

\put(23,5){\vector(1,0){10}}

\put(27,6){\shortstack{\mbox{$\sigma $}}}

\put(45,5){\makebox(0,0){\mbox{$ [l]$}}}

\put(75,5){\makebox(0,0){\mbox{$[k] $}}}

\put(57,5){\vector(1,0){10}}
\put(80,4){\makebox(0,0){\mbox{$.$}}}

\put(60,6){\shortstack{\mbox{$\omega $}}}

\put(57,8){\shortstack{\mbox{\small $ $}}}

\end{picture}}

 Let $\Sm $ be the subcategory of bijections in  $\Omega^s .$ This is a strict monoidal groupoid with tensor product $\oplus$ given by ordinal sum and with $[-1]$ as its unital object.

  A {\it right symmetric collection}  in a symmetric monoidal category $V$ is a functor $A:\Sm^{op}\rightarrow V .$ The value of $A$ on an object $[n]$ will be denoted $A_{n} .$ Notice, that this is not a standard operadic notation. Classically, the notation for $A_{[n]}$ is $A_{n+1}$ to stress the fact that $A_{n+1}$ is the space of operations of arity $n+1 .$  

The following definition is classical May definition \cite{May} of symmetric operad.
\begin{defin}\label{defsymop} A (right) symmetric  operad in $V$ is a right symmetric collection $A$ equipped with the following additional structure:

- a morphism  $e:I\rightarrow A_0$

- for every order preserving  map $\sigma:[n]\rightarrow [k]$ in $\Omega^s$   a morphism
: 
$$\mu_{\sigma}: A_{k}\otimes(A_{n_0}\otimes ... \otimes
A_{n_k})\longrightarrow A_{n},
$$
where $[n_i] = \sigma^{-1}(i).$ 

They must satisfy the following identities:
\begin{enumerate}\item for any composite of order preserving morphisms in $\Omega^s$ $$[n]\stackrel{\sigma}{\longrightarrow} [l] \stackrel{\omega}{\longrightarrow} [k] ,$$
the following diagram commutes

{\unitlength=1mm

\begin{picture}(300,45)(2,0)

\put(20,35){\makebox(0,0){\mbox{$\scriptstyle A_{k}\otimes
A_{l_{\bullet}}\otimes A_{n_0^{\bullet}} \otimes  ...
\otimes 
 A_{n_i^{\bullet}}\otimes  ... \otimes A_{n_k^{\bullet}}   
$}}}
\put(20,31){\vector(0,-1){12}}

\put(94,31){\vector(0,-1){12}}

\put(88,35){\makebox(0,0){\mbox{$\scriptstyle A_{k}\otimes
A_{l_{1}}\otimes A_{n_0^{\bullet}} \otimes  ...
\otimes A_{l_{i}}\otimes
 A_{n_i^{\bullet}}\otimes  ... \otimes A_{l_{k}}\otimes
A_{n_k^{\bullet}}   
$ }}}

\put(49,37){\makebox(0,0){\mbox{$\scriptstyle \simeq $}}}
\put(45,35.5){\vector(1,0){9}}

\put(20,15){\makebox(0,0){\mbox{$\scriptstyle A_{l}\otimes 
A_{n_0^{\bullet}} \otimes  ...
\otimes 
 A_{n_i^{\bullet}}\otimes  ... \otimes A_{n_k^{\bullet}}
$}}}

\put(94,15){\makebox(0,0){\mbox{$\scriptstyle A_{k}\otimes 
A_{n_{\bullet}} 
$}}}

\put(60,5){\makebox(0,0){\mbox{$ \scriptstyle A_{n} 
$}}}

\put(35,11){\vector(4,-1){19}}

\put(85,11){\vector(-4,-1){19}}

\end{picture}}

\noindent 
Here $$A_{l_{\bullet}}= A_{l_0}\otimes ...
\otimes A_{l_k},$$  
$$A_{n_{i}^{\bullet}} = A_{n_i^0} \otimes ...\otimes A_{n_i^{m_i}}$$
and $$ A_{n_{\bullet} } =  A_{n_0}\otimes ...
\otimes A_{n_k};$$

\item for an identity $\sigma = id : [n]\rightarrow [n]$ the diagram

{\unitlength=1mm
\begin{picture}(50,25)(30,2)

\put(97,20){\vector(-1,0){20}}

\put(60,17){\vector(0,-1){8}}

\put(60,20){\makebox(0,0){\mbox{\small$A_{n}\otimes 
A_{0}\otimes ... \otimes A_{0} 
$}}}

\put(114,20){\makebox(0,0){\mbox{\small$A_{n}\otimes 
{I}\otimes ... \otimes {I} 
$}}}

\put(60,5){\makebox(0,0){\mbox{\small$A_{n} 
$}}}

\put(105,15){\vector(-4,-1){30}}

\put(90,9){\makebox(0,0){\mbox{\small$id
$}}}

\end{picture}}

\noindent commutes;

\item for the unique morphism $[n]\rightarrow [0]$ the diagram

{\unitlength=1mm
\begin{picture}(50,25)(30,2)

\put(87,20){\vector(-1,0){15}}

\put(60,17){\vector(0,-1){8}}

\put(60,20){\makebox(0,0){\mbox{\small$A_{0}\otimes 
A_{n}
$}}}

\put(98,20){\makebox(0,0){\mbox{\small$I \otimes
A_{n}
$}}}

\put(60,5){\makebox(0,0){\mbox{\small$A_{n} 
$}}}

\put(95,17){\vector(-3,-1){25}}

\put(84,11){\makebox(0,0){\mbox{\small$id
$}}}

\end{picture}}

\noindent commutes.
\end{enumerate}
The following equivariance conditions are also required:
\begin{enumerate}      
 
\item For any order preserving $\sigma:[n]\rightarrow [k]$ and any bijection $\rho:[k]\rightarrow [k]$ the following diagram commutes:

{\unitlength=1mm
\begin{picture}(200,33)(-30,0)
\put(9,25){\makebox(0,0){\mbox{$A_{k}\otimes(A_{n_{\rho(0)}}
\otimes ... \otimes
A_{n_{\rho(k)}}) $}}}
\put(10,10){\vector(0,1){10}}
\put(-5,15){\shortstack{\mbox{$\scriptstyle A(\rho)\otimes\tau(\rho) $}}}
\put(32,25){\vector(1,0){12}}
\put(35,27){\shortstack{\mbox{$\mu_{\sigma} $}}}
\put(50,25){\makebox(0,0){\mbox{$ A_{n}$}}}
\put(50,10){\vector(0,1){10}}
\put(10,5){\makebox(0,0){\mbox{$A_{k}\otimes(A_{n_0}
\otimes ... \otimes
A_{n_k})$}}}
\put(31,5){\vector(1,0){13}}
\put(35,7){\shortstack{\mbox{$\mu_{\sigma} $}}}
\put(50,5){\makebox(0,0){\mbox{$A_{n}$}}}
\put(55,15){\makebox(0,0){\mbox{$\scriptstyle A(\pi) $}}}
\put(60,4){\shortstack{\mbox{$, $}}}
\end{picture}}
\noindent where $\tau(\rho)$ is the symmetry in $V$  which corresponds to permutation $\rho $ and $\pi = \Gamma_S(\rho;1,\ldots,1)$ is the permutation
, which permutes the fibers $[n_0],\ldots, [n_k]$ according to $\rho$ and whose restriction on each fiber is an identity.

\item For any order preserving $\sigma:[n]\rightarrow [k]$ and any set of bijections $\rho_i:[n_i]\rightarrow [n_i], 0\le i\le k ,$ 
the following diagram commutes 

{\unitlength=1mm
\begin{picture}(200,33)(-30,0)
\put(9,25){\makebox(0,0){\mbox{$A_{k}\otimes(A_{n_0}
\otimes ... \otimes
A_{n_{k}}) $}}}
\put(10,10){\vector(0,1){10}}
\put(-20,15){\shortstack{\mbox{$\scriptstyle 
id\otimes A(\rho_0)\otimes \ldots \otimes A(\rho_k) $}}}
\put(32,25){\vector(1,0){12}}
\put(35,27){\shortstack{\mbox{$\mu_{\sigma} $}}}
\put(50,25){\makebox(0,0){\mbox{$ A_{n}$}}}
\put(50,10){\vector(0,1){10}}
\put(10,5){\makebox(0,0){\mbox{$A_{k}\otimes(A_{n_0}
\otimes ... \otimes
A_{n_k})$}}}
\put(31,5){\vector(1,0){13}}
\put(35,7){\shortstack{\mbox{$\mu_{\sigma} $}}}
\put(50,5){\makebox(0,0){\mbox{$A_{n}$}}}
\put(62,15){\makebox(0,0){\mbox{$\scriptstyle A(\rho_0\oplus\ldots\oplus \rho_k) $}}}
\put(60,4){\shortstack{\mbox{$, $}}}
\end{picture}}
\noindent 
\end{enumerate}

\end{defin}

We can give an alternative definition of  symmetric operad
\cite{EHBat}.

\begin{defin}\label{symop}

 A (right) symmetric  operad in $V$ is a right symmetric collection $A$ equipped with the following additional structure:

- a morphism  $e:I\rightarrow A_0$

- for every order preserving  map $\sigma:[n]\rightarrow [k]$ in $\Omega^s$   a morphism:

$$\mu_{\sigma}: A_{k}\otimes(A_{n_0}\otimes ... \otimes
A_{n_k})\longrightarrow A_{n},
$$
where $[n_i] = \sigma^{-1}(i).$ 

They must satisfy the same conditions as in the definition \ref{defsymop} with respect to order preserving maps and identities but the equivariance conditions are replaced by the following:
\begin{enumerate}

\item

 For every commutative diagram in $\Omega^s $

{\unitlength=1mm

\begin{picture}(40,28)(-29,2)

\put(13,25){\makebox(0,0){\mbox{$ [n']$}}}
\put(13,21){\vector(0,-1){10}}
\put(10,15){\shortstack{\mbox{$\pi $}}}
\put(42,15){\shortstack{\mbox{$\rho $}}}

\put(22,25){\vector(1,0){10}}

\put(26,26){\shortstack{\mbox{$\sigma' $}}}

\put(41,25){\makebox(0,0){\mbox{$ [k']$}}}
\put(41,21){\vector(0,-1){10}}

\put(13,7){\makebox(0,0){\mbox{$ [n]$}}}

\put(23,7){\vector(1,0){10}}

\put(27,8){\shortstack{\mbox{$\sigma $}}}

\put(41,7){\makebox(0,0){\mbox{$[k]$}}}

\end{picture}}

whose vertical maps are bijections and whose horizontal
maps are order preserving the following diagram commutes:

{\unitlength=1mm

\begin{picture}(200,33)(-30,0)

\put(9,25){\makebox(0,0){\mbox{$A_{k'}\otimes(A_{n'_{\rho(0)}}
\otimes ... \otimes
A_{n'_{\rho(k)}}) $}}}
\put(10,10){\vector(0,1){10}}
\put(-5,15){\shortstack{\mbox{$\scriptstyle A(\rho)\otimes\tau(\rho) $}}}

\put(32,25){\vector(1,0){12}}

\put(35,27){\shortstack{\mbox{$\mu_{\sigma'} $}}}

\put(50,25){\makebox(0,0){\mbox{$ A_{n'}$}}}
\put(50,10){\vector(0,1){10}}

\put(10,5){\makebox(0,0){\mbox{$A_{k}\otimes(A_{n_0}
\otimes ... \otimes
A_{n_k})$}}}

\put(31,5){\vector(1,0){13}}

\put(35,7){\shortstack{\mbox{$\mu_{\sigma} $}}}

\put(50,5){\makebox(0,0){\mbox{$A_{n}$}}}

\put(55,15){\makebox(0,0){\mbox{$\scriptstyle A(\pi) $}}}

\put(60,4){\shortstack{\mbox{$, $}}}

\end{picture}}
\noindent where $\tau(\rho)$ is the symmetry in $V$  which corresponds to permutation $\rho .$ 

\item For every commutative diagram in $\Omega^s$ 

{\unitlength=1mm

\begin{picture}(40,28)(-29,2)

\put(13,25){\makebox(0,0){\mbox{$ [n'']$}}}
\put(13,21){\vector(0,-1){10}}
\put(10,15){\shortstack{\mbox{$\sigma $}}}
\put(42,15){\shortstack{\mbox{$\eta' $}}}

\put(22,25){\vector(1,0){10}}

\put(26,26){\shortstack{\mbox{$\sigma' $}}}

\put(41,25){\makebox(0,0){\mbox{$ [n']$}}}
\put(41,21){\vector(0,-1){10}}

\put(13,7){\makebox(0,0){\mbox{$ [n]$}}}

\put(23,7){\vector(1,0){10}}

\put(27,8.5){\shortstack{\mbox{$\eta $}}}

\put(41,7){\makebox(0,0){\mbox{$[k]$}}}

\end{picture}}

\noindent where $\sigma,\sigma'$ are bijections and $\eta,\eta'$
are order preserving maps, the following diagram commutes 

{\unitlength=1mm

\begin{picture}(200,53)(-30,0)

\put(10,25){\makebox(0,0){\mbox{$A_{k}\otimes (A_{n''_0}
\otimes ... \otimes
A_{n''_k}) $}}}
\put(10,8){\vector(0,1){12}}
\put(-18,13){\shortstack{\mbox{$\scriptstyle 1\otimes A(\sigma_0)\otimes\ldots\otimes A(\sigma_k) $}}}

\put(-18,34){\shortstack{\mbox{$\scriptstyle 1\otimes A(\sigma'_0)\otimes\ldots\otimes A(\sigma'_k) $}}}

\put(46,34){\shortstack{\mbox{$\scriptstyle  A(\sigma') $}}}
\put(46,13){\shortstack{\mbox{$\scriptstyle  A(\sigma) $}}}
\put(35,47){\makebox(0,0){\mbox{$ \mu_{\eta'}$}}}
\put(34,7){\makebox(0,0){\mbox{$ \mu_{\eta}$}}}

\put(45,25){\makebox(0,0){\mbox{$ A_{n''}$}}}
\put(45,8){\vector(0,1){12}}

\put(75,25){\makebox(0,0){\mbox{$ $}}}

\put(57,21){\shortstack{\mbox{$ $}}}

\put(10,5){\makebox(0,0){\mbox{$A_{k}\otimes (A_{n_0}
\otimes ... \otimes
A_{n_k})$}}}

\put(29,5){\vector(1,0){10}}

\put(27,6){\shortstack{\mbox{$ $}}}

\put(45,5){\makebox(0,0){\mbox{$A_{n}$}}}

\put(75,5){\makebox(0,0){\mbox{$ $}}}


\put(60,6){\shortstack{\mbox{$ $}}}


\put(10,45){\makebox(0,0){\mbox{$A_{k}\otimes(A_{n'_0}
\otimes ... \otimes
A_{n'_k}) $}}}

\put(45,45){\makebox(0,0){\mbox{$ A_{n'}$}}}

\put(29,45){\vector(1,0){10}}

\put(10,41){\vector(0,-1){12}}
\put(45,41){\vector(0,-1){12}}

\end{picture}}

\end{enumerate}

\end{defin}

\begin{pro} The definition \ref{defsymop} and \ref{symop}  are equivalent.
 \end{pro}

We leave this proposition as an exercise for the reader. 

\

Let $\Br$ be the groupoid of braid groups. We will regard the objects of $\Br$ as ordinals. There is a monoidal structure on $\Br$ given by ordinal sum on objects and concatenation of braids on morphism. The ordinal $[-1]$ is the unital object.

The following is the definition of braided operad from \cite{fied}.
 A {\it right braided collection}  in a symmetric monoidal category $V$ is a functor $A:\Br^{op}\rightarrow V .$ The value of $A$ on an object $[n]$ will be denoted $A_n .$ 
\begin{defin}\label{defbrop} A right braided  operad in $V$ is a right braided collection $A$ equipped with the following additional structure:

- a morphism  $e:I\rightarrow A_0$

- for every order preserving  map $\sigma:[n]\rightarrow [k]$ in $\Omega^s$   a morphism
: 
$$\mu_{\sigma}: A_{k}\otimes(A_{n_0}\otimes ... \otimes
A_{n_k})\longrightarrow A_{n},
$$
where $[n_i] = \sigma^{-1}(i).$ 

They must satisfy the   identities (1-3) from the definition \ref{defsymop}
and  the following 
  two equivariancy conditions:

\begin{enumerate}      
 
\item For any order preserving $\sigma:[n]\rightarrow [k]$ and any braid $\rho:[k]\rightarrow [k]$ the following diagram commutes:

{\unitlength=1mm
\begin{picture}(200,33)(-30,0)
\put(9,25){\makebox(0,0){\mbox{$A_{k}\otimes(A_{n_{\rho(0)}}
\otimes ... \otimes
A_{n_{\rho(k)}}) $}}}
\put(10,10){\vector(0,1){10}}
\put(-5,15){\shortstack{\mbox{$\scriptstyle A(\rho)\otimes\tau(\rho) $}}}
\put(32,25){\vector(1,0){12}}
\put(35,27){\shortstack{\mbox{$\mu_{\sigma} $}}}
\put(50,25){\makebox(0,0){\mbox{$ A_{n}$}}}
\put(50,10){\vector(0,1){10}}
\put(10,5){\makebox(0,0){\mbox{$A_{k}\otimes(A_{n_0}
\otimes ... \otimes
A_{n_k})$}}}
\put(31,5){\vector(1,0){13}}
\put(35,7){\shortstack{\mbox{$\mu_{\sigma} $}}}
\put(50,5){\makebox(0,0){\mbox{$A_{n}$}}}
\put(55,15){\makebox(0,0){\mbox{$\scriptstyle A(\pi) $}}}
\put(60,4){\shortstack{\mbox{$, $}}}
\end{picture}}
\noindent where $\tau(\rho)$ is the symmetry in $V$  which corresponds to the braid $\rho $ and $\pi = \Gamma_B(\rho;1,\ldots,1)$ is 
a braid obtained from $\rho$ by replacing the $i$-th strand of $\rho$ by $n_i$ parallel strands for each $i .$

\item For any order preserving $\sigma:[n]\rightarrow [k]$ and any set of braids $\rho_i:[n_i]\rightarrow [n_i], 0\le i\le k ,$ 
the following diagram commutes 

{\unitlength=1mm
\begin{picture}(200,33)(-30,0)
\put(9,25){\makebox(0,0){\mbox{$A_{k}\otimes(A_{n_0}
\otimes ... \otimes
A_{n_{k}}) $}}}
\put(10,10){\vector(0,1){10}}
\put(-20,15){\shortstack{\mbox{$\scriptstyle 
id\otimes A(\rho_0)\otimes \ldots \otimes A(\rho_k) $}}}
\put(32,25){\vector(1,0){12}}
\put(35,27){\shortstack{\mbox{$\mu_{\sigma} $}}}
\put(50,25){\makebox(0,0){\mbox{$ A_{n}$}}}
\put(50,10){\vector(0,1){10}}
\put(10,5){\makebox(0,0){\mbox{$A_{k}\otimes(A_{n_0}
\otimes ... \otimes
A_{n_k})$}}}
\put(31,5){\vector(1,0){13}}
\put(35,7){\shortstack{\mbox{$\mu_{\sigma} $}}}
\put(50,5){\makebox(0,0){\mbox{$A_{n}$}}}
\put(62,15){\makebox(0,0){\mbox{$\scriptstyle A(\rho_0\oplus\ldots\oplus \rho_k) $}}}
\put(60,4){\shortstack{\mbox{$, $}}}
\end{picture}}
\noindent 
\end{enumerate}

\end{defin}

\section{ $n$-ordinals and  quasibijections} \label{nord}

 \begin{defin} An $n$-ordinal consists of a finite set $T$  equipped with $n$ binary relations  $<_0, \ldots, <_{n-1} $
  satisfying the following axioms
  \begin{enumerate}
  \item  \ $<_p$ is nonreflexive;
 \item \ for every pair $a,b$ of distinct elements of $T$ there exists exactly one $p$ such that
 $$a<_p b \ \ \mbox{or} \ \ b<_p a ;$$
 \item \  if $a<_p b $ \ and \ $b<_q c$ \ then\
 $a<_{min(p,q)} c .$
 \end{enumerate}
 \end{defin}
 
 Every $n$-ordinal can be represented as a pruned planar tree with $n$ levels.   For example, the $2$-ordinal
 \begin{equation}\label{ordinal} 0<_0 1,\ 0<_0 2, \ 0<_0 3,\ 
1<_1 2,\ 2<_1 3,\ 2<_1 3\end{equation} 
is represented by the following pruned tree

{\unitlength=0.71mm
\begin{picture}(100,20)(-15,-1)

\put(60,0){\line(-1,1){5}}
\put(60,0){\line(1,1){5}}
\put(55,5){\line(0,1){5}} \put(55,13){\makebox(0,0){\mbox{$0$}}}
\put(65,5){\line(1,1){5}}\put(60,13){\makebox(0,0){\mbox{$1$}}}
\put(65,5){\line(-1,1){5}}\put(65,13){\makebox(0,0){\mbox{$2$}}}
\put(65,5){\line(0,1){5}}\put(70,13){\makebox(0,0){\mbox{$3$}}}

\end{picture}}

\noindent
 See \cite{SymBat} for a more detailed discussion.  
 \begin{defin}\label{mapofordinals}  \  A map of $n$-ordinals
 $$\sigma: T \rightarrow S$$   is a map $\sigma:T\rightarrow S$ of underlying sets such that  $$i<_p j \ \mbox{in} \ T $$   implies that
 \begin{enumerate}
\item \ $\sigma(i) <_r \sigma(j)$   for some $r\ge p$ or
\item \  $\sigma(i)= \sigma(j)$ or
\item $\sigma(j) <_r \sigma(i)$ for $r>p .$
\end{enumerate}
 \end{defin}

For every $i\in S$ the preimage $\sigma^{-1}(i)$ ({\it the fiber of $\sigma$ over $i$}) has a natural structure of an $n$-ordinal.

We denote by $Ord(n)$ {\it the skeletal category of $n$-ordinals .} The category $Ord(n)$   is monoidal. The monoidal structure $\oplus$ is defined as follows. For two $n$-ordinals
$S$ and $T$ the $n$-ordinal $S\oplus T$ has as an underlying set the union of underlying sets of $S$ and $T .$ The orders $<_k$ restricted to the elements of $S$ and $T$ coincide with respective orders on $S$ and $T .$ and $a<_0 b$ if $a\in S$ and $b\in T .$ The unital object for this monoidal structure is empty $n$-ordinal.

An $n$-ordinal structure on $T$ determines a linear order  (called {\it total order}) on the elements of $T$ 
 as follows:
$$a<b \ \ \ \mbox{iff} \ \ \ a<_r b \ \ \ \mbox{for some} \ \ 0\le r \le n-1 \ .$$
We will denote by $[T]$ the set $T$ with its total linear order.
In this way we have a monoidal functor $$[-]:Ord(n)\rightarrow \Omega^s.$$  This functor is faithful but not full. For example,  no morphism from the $2$-ordinal (\ref{ordinal}) to the $2$-ordinal $0<_1 1$ can reverse the order of  $1,2$ and $3$ 

We also introduce the category of $\infty$-ordinals $Ord(\infty).$ 

 \begin{defin} An $\infty$-ordinal consists of a finite set $T$  equipped with a sequence of binary relations  $<_0, <_{-1}, <_{-2},\ldots $
  satisfying the following axioms
  \begin{enumerate}
  \item  \ $<_p$ is nonreflexive;
 \item \ 
for every pair $a,b$ of distinct elements of $T$ there exists exactly one $p$ such that
 $$a<_p b \ \ \mbox{or} \ \ b<_p a ;$$
 \item \  if $a<_p b $ \ and \ $b<_q c$ \ then\
 $a<_{min(p,q)} c .$
 \end{enumerate}
 \end{defin}

The definition of morphism between $\infty$-ordinals coincides with the Definition \ref{mapofordinals}. 
The category $Ord(\infty)$ is {\it the skeletal category of $\infty$-ordinals .} As for $Ord(n)$ we have a functor of total order $$[-]: Ord(\infty)\rightarrow \Omega^s .$$

For a $k$-ordinal $R\ ,\ k\le n $ we consider its $(n-k)$-th {\it vertical suspension} $S^{n-k}R$ which is an $n$-ordinal with the underlying set $R ,$ and the order $<_m $ equal the  order $<_{m-k}$ on $R$   (so   $<_m$ are empty for $0\le m <n-k .$)
We also can consider {\it 
the horizontal $(n-1)$-suspension} $T^{n-k}R$  which is a $n$-ordinal with the underlying set $R ,$ and the order $<_m$ equal the order on $R$ (so   $<_m$ are empty for $k-1< m \le n-1 .$)  

The  vertical suspension provides us with a functor
$S:Ord(n)\rightarrow Ord(n+1) .$ We also define an $\infty$-suspension functor $Ord(n)\rightarrow Ord(\infty)$ as follows. For an $n$-ordinal $T$ its $\infty$-suspension is an $\infty$-ordinal $S^{\infty}T$ whose underlying set is the same as the underlying set of $T$ and $a<_p b$ in $S^{\infty}T$ if 
$a<_{n+p-1} b$ in $T .$ It is not hard to see that the sequence 
$$ Ord(0)\stackrel{S}{\longrightarrow} Ord(1) \stackrel{S}{\longrightarrow} Ord(2) \longrightarrow \ldots \stackrel{S}{\longrightarrow} Ord(n) \longrightarrow \ldots \stackrel{S^{\infty}}{\longrightarrow} Ord(\infty),$$
 exhibits $Ord(\infty)$ as a colimit of $Ord(n) .$

\begin{defin} A map of  $n$-ordinals is called a quasibijection if it is a bijection of the underlying sets.\end{defin} 

Let $Q_n \ , 1\le n\le\infty $ be {\it the subcategory of quasibijections} of $Ord(n) .$ The total order functor induces then a functor which we will denote by the same symbol: 
$$[-]:Q_n\rightarrow \Sm .$$ 

\begin{defin}A map $\sigma$ of $n$-ordinals $1\le n\le \infty$ is called order preserving  if it preserves the total  orders  in the usual sense or equivalently  only  conditions $1$ and $2$ from the Definition \ref{mapofordinals} hold for $\sigma.$ \end{defin}

\begin{lem}\label{factorisation} For every morphism $\sigma:T\rightarrow S$ in $Ord(n)$ $1\le n\le \infty$ there exists a
 factorisation
$$T\stackrel{\pi}{\longrightarrow} T' \stackrel{\nu}{\longrightarrow} S$$
where $\pi$ is a quasibijection, $\nu$ is order preserving and $\pi$ preserves total order on fibers 
of $\nu$.
\end{lem}

\Proof For $n=1$ this factorisation is trivial, since all maps of $1$-ordinals are order preserving.

Let $n=2 .$ Let $\sigma:T\rightarrow S$ be a map of $2$-ordinals and let $S= S[k]$ be a suspension of the $1$-ordinal $[k] .$ Let
$T'$ be the $2$-ordinal whose underlying set is the same as that of $T,$  whose only nonempty order is $<_1 $ and whose total order coincides with $[T] .$  So $T'$ itself is a vertically suspended $1$-ordinal. Now, one can factorise the map
$[\sigma]:[T]\rightarrow [S]$ in $\Omega^s $ 
 $$[T]\stackrel{\pi}{\longrightarrow} [T'] \stackrel{\nu}{\longrightarrow} [S] $$
with $\nu$ being total order preserving and $\pi$ a bijection
which preserves the order on the fibers of $\sigma $ \cite{EHBat}. Obviously, $\nu$ can be considered as a map of $2$-ordinals and it is order preserving.
Let us check that $\pi$ is also a map of $2$-ordinals.  
Indeed, if $i,j$ are from the same fiber of $\sigma$ then 
$\pi$ preserves their order. If $i<_0 j$  in $T$ and they are from different 
fibers then there is no restriction on $\pi$ since $T'$ is a suspended $1$-ordinal. Finally, if $i<_1 j$  in $T$ and they are from different 
fibers then $\sigma(i)<_1\sigma(j), $ so $\pi(i)<_1\pi(j) $
 because $\nu$ is order preserving.

Finally, if $S$ is an arbitrary $2$-ordinal then 
$S = S_1\oplus\ldots \oplus S_k$ for some suspended $1$-ordinals $S_1,\ldots,S_k$ and moreover,     
$$\sigma = \sigma_1\oplus\ldots \oplus \sigma_k
:T= T_1\oplus\ldots \oplus T_k \rightarrow  S_1\oplus\ldots \oplus S_k .$$ 
By applying the previous result to each $\sigma_k$ we obtain a required factorisation of $\sigma .$ 

The factorisation for $n>2$ can be obtained similarly. 

\Q

\section{Quasisymmetric $n$-operads.}\label{qoperad}

We now recall  the definition of pruned $(n-1)$-terminal  $n$-operad \cite{SymBat}. Since we do not need other types of
$n$-operads in this paper we will call them simply $n$-operads. 
The notation $U_n$ means the terminal $n$-ordinal.  

Let $V$ be a symmetric monoidal category. For a morphism of $n$-ordinals $\sigma:T\rightarrow S$ the $n$-ordinal $T_i$ is the fiber $\sigma^{-1}(i) .$ 
\begin{defin}\label{defnoper}  An $n$-operad in $V$ is 
a collection $A_T, \ T\in Ord(n)$ of objects of $V$ equipped with the following structure :

- a morphism $e: I \rightarrow  A_{U_n}$ (the unit);

- for every morphism $\sigma:T \rightarrow S$ in $Ord(n) ,$ 
a morphism 
$$m_{\sigma}:A_S\otimes A_{T_0}\otimes ... \otimes A_{T_k}
 \rightarrow A_T \mbox{\ \ (the multiplication}).$$

They must satisfy the following identities:

- for any composite $T\stackrel{\sigma}{\rightarrow} S \stackrel{\omega}{\rightarrow} R ,$
the associativity diagram

{\unitlength=1mm

\begin{picture}(300,45)(2,0)

\put(20,35){\makebox(0,0){\mbox{$\scriptstyle A_R\otimes
A_{S_{\bullet}}\otimes A_{T_0^{\bullet}} \otimes  ...
\otimes 
 A_{T_i^{\bullet}}\otimes  ... \otimes A_{T_k^{\bullet}}   
$}}}
\put(20,31){\vector(0,-1){12}}

\put(94,31){\vector(0,-1){12}}

\put(88,35){\makebox(0,0){\mbox{$\scriptstyle A_R\otimes
A_{S_{0}}\otimes A_{T_1^{\bullet}} \otimes  ...
\otimes A_{S_{i}}\otimes
 A_{T_i^{\bullet}}\otimes  ... \otimes A_{S_{k}}\otimes
A_{T_k^{\bullet}}   
$ }}}

\put(50,35){\makebox(0,0){\mbox{$\scriptstyle \simeq $}}}

\put(20,15){\makebox(0,0){\mbox{$\scriptstyle A_S\otimes 
A_{T_1^{\bullet}} \otimes  ...
\otimes 
 A_{T_i^{\bullet}}\otimes  ... \otimes A_{T_k^{\bullet}}
$}}}

\put(94,15){\makebox(0,0){\mbox{$\scriptstyle A_R\otimes 
A_{T_{\bullet}} 
$}}}

\put(60,5){\makebox(0,0){\mbox{$ \scriptstyle A_T 
$}}}

\put(35,11){\vector(4,-1){19}}

\put(85,11){\vector(-4,-1){19}}

\end{picture}}

\noindent commutes,
where $$A_{S_{\bullet}}= A_{S_0}\otimes ...
\otimes A_{S_k},$$  
$$A_{T_{i}^{\bullet}} = A_{T_i^0} \otimes ...\otimes A_{T_i^{m_i}}$$
and $$ A_{T_{\bullet} } =  A_{T_0}\otimes ...
\otimes A_{T_k};$$

- for an identity $\sigma = id : T\rightarrow T$ the diagram

{\unitlength=1mm
\begin{picture}(50,25)(30,2)

\put(97,20){\vector(-1,0){20}}

\put(60,17){\vector(0,-1){8}}

\put(60,20){\makebox(0,0){\mbox{\small$A_T\otimes 
A_{U_n}\otimes ... \otimes A_{U_n} 
$}}}

\put(114,20){\makebox(0,0){\mbox{\small$A_T\otimes 
{I}\otimes ... \otimes {I} 
$}}}

\put(60,5){\makebox(0,0){\mbox{\small$A_T 
$}}}

\put(105,15){\vector(-4,-1){30}}

\put(90,9){\makebox(0,0){\mbox{\small$id
$}}}

\end{picture}}

\noindent commutes;

- for the unique morphism $T\rightarrow U_n$ the diagram

{\unitlength=1mm
\begin{picture}(50,25)(30,2)

\put(87,20){\vector(-1,0){15}}

\put(60,17){\vector(0,-1){8}}

\put(60,20){\makebox(0,0){\mbox{\small$A_{U_n}\otimes 
A_T
$}}}

\put(98,20){\makebox(0,0){\mbox{\small$I \otimes
A_T
$}}}

\put(60,5){\makebox(0,0){\mbox{\small$A_T 
$}}}

\put(95,17){\vector(-3,-1){25}}

\put(84,11){\makebox(0,0){\mbox{\small$id
$}}}

\end{picture}}

\noindent commutes.

\end{defin}

Let $\sigma:T\rightarrow S$ be a quasibijection and $A$ be a pruned $n$-operad.
Since a fiber of $\sigma$ is the terminal $n$-ordinal $U_n ,$  the multiplication
$$\mu_{\sigma}: A_{S}\otimes(A_{U_n}\otimes ... \otimes
A_{U_n})\longrightarrow A_{T}$$
in composition with the morphism
$$ A_S\rightarrow  A_{S}\otimes(I\otimes ... \otimes
I) \rightarrow  A_{S}\otimes(A_{U_n}\otimes ... \otimes
A_{U_n})$$ induces a morphism
$$ A(\sigma): A_S\rightarrow A_T.$$
It is not hard to see  that in this way $A$ becomes a contravariant  functor on  $Q_n .$  

\begin{defin} We call a pruned $n$-operad $A$ {\it quasisymmetric} if  for every quasibijection  $\sigma:T\rightarrow S$  the  morphism
$$ A(\sigma): A_S\rightarrow A_T$$
is an isomorphism.

\end{defin}

The desymmetrisation functor from symmetric to $n$-operads for finite $n$ was defined in \cite{EHBat} using pulling back along the functor $[-]:Ord(n)\rightarrow \Omega^s .$  It was shown that this functor has a left adjoint which we call symmetrisation. We can obviously extend these definitions to $n=\infty.$  By construction the desymmetrisation of a symmetric operad is a quasisymmetric $n$-operad for any $n.$

 Let $\Pi Q_n$ be the fundamental groupoid of $Q_n$.   A quasisymmetric operad provides, therefore , a contravariant functor on $\Pi Q_n$.

\begin{defin} A $Q_n$-collection is a contravariant functor on $ Q_n$. A $\Pi Q_n$-collection is a contravariant functor on $\Pi Q_n .$\end{defin}

\begin{defin} A $Q_n$-operad is a $\Pi Q_n$-collection $A$ together with
the following structure
\begin{itemize}
 \item for every order preserving  map $\sigma:
T\rightarrow S$ the  usual  operadic map: 
$$\mu_{\sigma}: A_{S}\otimes (A_{T_0}\otimes ... \otimes
A_{T_k})\longrightarrow A_{T}.$$
This collection of maps must satisfy the usual associativity and unitarity conditions plus
 two equivariancy conditions:
\item For every commutative diagram

{\unitlength=0.7mm

\begin{picture}(58,35)(-35,-2)

\put(10,25){\makebox(0,0){\mbox{$ T'$}}}
\put(10,20){\vector(0,-1){10}}
\put(12,15){\shortstack{\mbox{$ $}}}

\put(22,25){\vector(1,0){10}}

\put(24,26){\shortstack{\mbox{$\sigma' $}}}

\put(45,25){\makebox(0,0){\mbox{$ S'$}}}
\put(45,20){\vector(0,-1){10}}

\put(75,25){\makebox(0,0){\mbox{$ $}}}


\put(57,21){\shortstack{\mbox{$ $}}}


\put(10,5){\makebox(0,0){\mbox{$ T$}}}

\put(23,5){\vector(1,0){10}}

\put(27,6){\shortstack{\mbox{$\sigma $}}}

\put(45,5){\makebox(0,0){\mbox{$S$}}}

\put(75,5){\makebox(0,0){\mbox{$ $}}}


\put(60,6){\shortstack{\mbox{$ $}}}


\end{picture}}

where vertical maps are quasibijections and horizontal
maps are  order preserving the diagram

{\unitlength=1mm

\begin{picture}(200,33)(-24,-1)

\put(10,25){\makebox(0,0){\mbox{$A_{S}\otimes(A_{T_0}
\otimes ... \otimes
A_{T_k}) $}}}
\put(10,20){\vector(0,-1){10}}
\put(12,15){\shortstack{\mbox{$ $}}}

\put(29,25){\vector(1,0){10}}

\put(24,26){\shortstack{\mbox{$ $}}}

\put(45,25){\makebox(0,0){\mbox{$ A_T$}}}
\put(45,20){\vector(0,-1){10}}

\put(75,25){\makebox(0,0){\mbox{$ $}}}


\put(57,21){\shortstack{\mbox{$ $}}}


\put(10,5){\makebox(0,0){\mbox{$A_{S'}\otimes(A_{T'_0}
\otimes ... \otimes
A_{T'_k})$}}}

\put(29,5){\vector(1,0){10}}

\put(27,6){\shortstack{\mbox{$ $}}}

\put(45,5){\makebox(0,0){\mbox{$A_{T'}$}}}

\put(75,5){\makebox(0,0){\mbox{$ $}}}


\put(60,6){\shortstack{\mbox{$ $}}}


\end{picture}}
 commutes

\item  For every commutative diagram

{\unitlength=0.7mm

\begin{picture}(200,35)(-42,-2)

\put(10,25){\makebox(0,0){\mbox{$ T$}}}
\put(10,20){\vector(0,-1){10}}
\put(12,15){\shortstack{\mbox{$ $}}}

\put(22,25){\vector(1,0){10}}

\put(24,26){\shortstack{\mbox{$\sigma' $}}}

\put(45,25){\makebox(0,0){\mbox{$ T'$}}}
\put(45,20){\vector(0,-1){10}}

\put(7,15){\makebox(0,0){\mbox{$\sigma $}}}


\put(57,21){\shortstack{\mbox{$ $}}}


\put(50,15){\makebox(0,0){\mbox{$\eta'$}}}

\put(10,5){\makebox(0,0){\mbox{$ T''$}}}

\put(23,5){\vector(1,0){10}}

\put(27,6){\shortstack{\mbox{$\eta $}}}

\put(45,5){\makebox(0,0){\mbox{$S$}}}

\put(75,5){\makebox(0,0){\mbox{$ $}}}


\put(60,6){\shortstack{\mbox{$ $}}}


\end{picture}}
where $\sigma,\sigma'$ are quasibijections and $\eta,\eta'$
are order preserving, the diagram  

{\unitlength=0.9mm

\begin{picture}(200,53)(-30,0)

\put(10,25){\makebox(0,0){\mbox{$A_{S}\otimes(A_{T_0}
\otimes ... \otimes
A_{T_k}) $}}}
\put(10,10){\vector(0,1){10}}
\put(12,15){\shortstack{\mbox{$ $}}}


\put(24,26){\shortstack{\mbox{$ $}}}

\put(45,25){\makebox(0,0){\mbox{$ A_T$}}}
\put(45,10){\vector(0,1){10}}

\put(75,25){\makebox(0,0){\mbox{$ $}}}


\put(57,21){\shortstack{\mbox{$ $}}}


\put(10,5){\makebox(0,0){\mbox{$A_{S}\otimes(A_{T''_0}
\otimes ... \otimes
A_{T''_k})$}}}

\put(31,5){\vector(1,0){10}}

\put(27,6){\shortstack{\mbox{$ $}}}

\put(45,5){\makebox(0,0){\mbox{$A_{T''}$}}}

\put(75,5){\makebox(0,0){\mbox{$ $}}}


\put(60,6){\shortstack{\mbox{$ $}}}


\put(10,45){\makebox(0,0){\mbox{$A_{S}\otimes(A_{T'_0}
\otimes ... \otimes
A_{T'_k}) $}}}

\put(45,45){\makebox(0,0){\mbox{$ A_{T'}$}}}

\put(31,45){\vector(1,0){10}}

\put(10,39){\vector(0,-1){10}}
\put(45,39){\vector(0,-1){10}}

\end{picture}}

\noindent  commutes. \end{itemize}
\end{defin}

\begin{theorem} The category of $Q_n$-operads is equivalent  to the category of quasisymmetric $n$-operads. \end{theorem}

 \Proof Obviously, every quasisymmetric $n$-operad is a $Q_n$-operad.   Let us construct an inverse   functor. Given a $Q_n$-operad  $C$
we define a quasisymmetric operad $A$ on an $n$-ordinal $T$ to be equal   to $C_T$. We have to define $A$ on an arbitrary  map of $n$-ordinals 
$\sigma:T\rightarrow S$.

 Let us choose a factorisation of $\sigma$ according to  Lemma \ref{factorisation}.

Now we can define operadic multiplication by the following commutative diagram

{\unitlength=1mm

\begin{picture}(200,33)(-25,-2)

\put(10,25){\makebox(0,0){\mbox{$A_{S}\otimes (A_{T_0}
\otimes ... \otimes
A_{T_k}) $}}}
\put(10,20){\vector(0,-1){10}}
\put(-15,15){\shortstack{\mbox{$\scriptstyle 1\otimes (\alpha^{-1}_{\pi_1}\otimes
\ldots\otimes \alpha^{-1}_{\pi_k})$}}}

\put(29,25){\vector(1,0){10}}

\put(32,26){\shortstack{\mbox{$\scriptstyle \mu_{\sigma} $}}}

\put(45,25){\makebox(0,0){\mbox{$ A_T$}}}
\put(45,10){\vector(0,1){10}}
\put(47,15){\shortstack{\mbox{$\scriptstyle \alpha_{\pi} $}}}

\put(75,25){\makebox(0,0){\mbox{$ $}}}


\put(57,21){\shortstack{\mbox{$ $}}}


\put(10,5){\makebox(0,0){\mbox{$ A_{S}\otimes(A_{T'_0}
\otimes ... \otimes
A_{T'_k})$}}}

\put(29,5){\vector(1,0){10}}

\put(32,6){\shortstack{\mbox{$\scriptstyle \mu_\nu $}}}

\put(45,5){\makebox(0,0){\mbox{$A_{T'}$}}}

\put(75,5){\makebox(0,0){\mbox{$ $}}}


\put(60,6){\shortstack{\mbox{$ $}}}


\end{picture}}

The second equivariancy axiom implies that this definition does not depend on a chosen factorisation.   
Suppose now we have a composite 
$$T\stackrel{\sigma}{\longrightarrow} S \stackrel{\omega}{\longrightarrow} R .$$
It generates the following factorization diagram

{\unitlength=0.8mm

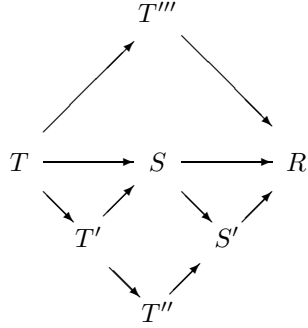
\begin{picture}(200,62)(-25,-2)

\put(20,25){\makebox(0,0){\mbox{$T$}}}

\put(24,25){\vector(1,0){15}}
\put(24,30){\vector(1,1){15}}
\put(47,46){\vector(1,-1){15}}

\put(24,20){\vector(1,-1){5}}
\put(47,20){\vector(1,-1){5}}
\put(35,7.5){\vector(1,-1){5}}

\put(45,3){\vector(1,1){5}}
\put(57,15){\vector(1,1){5}}
\put(34,16){\vector(1,1){5}}

\put(43,25){\makebox(0,0){\mbox{$S$}}}
\put(43,50){\makebox(0,0){\mbox{$T'''$}}}
\put(43,0){\makebox(0,0){\mbox{$T''$}}}

\put(31.5,12.5){\makebox(0,0){\mbox{$T'$}}}
\put(54.5,12.5){\makebox(0,0){\mbox{$S'$}}}

\put(47,25){\vector(1,0){15}}

\put(66,25){\makebox(0,0){\mbox{$R$}}}
\put(57,21){\shortstack{\mbox{$ $}}}

\end{picture}}

\

\noindent which in its turn generates the following huge diagram
\vspace{10mm}
{\unitlength=0.9mm

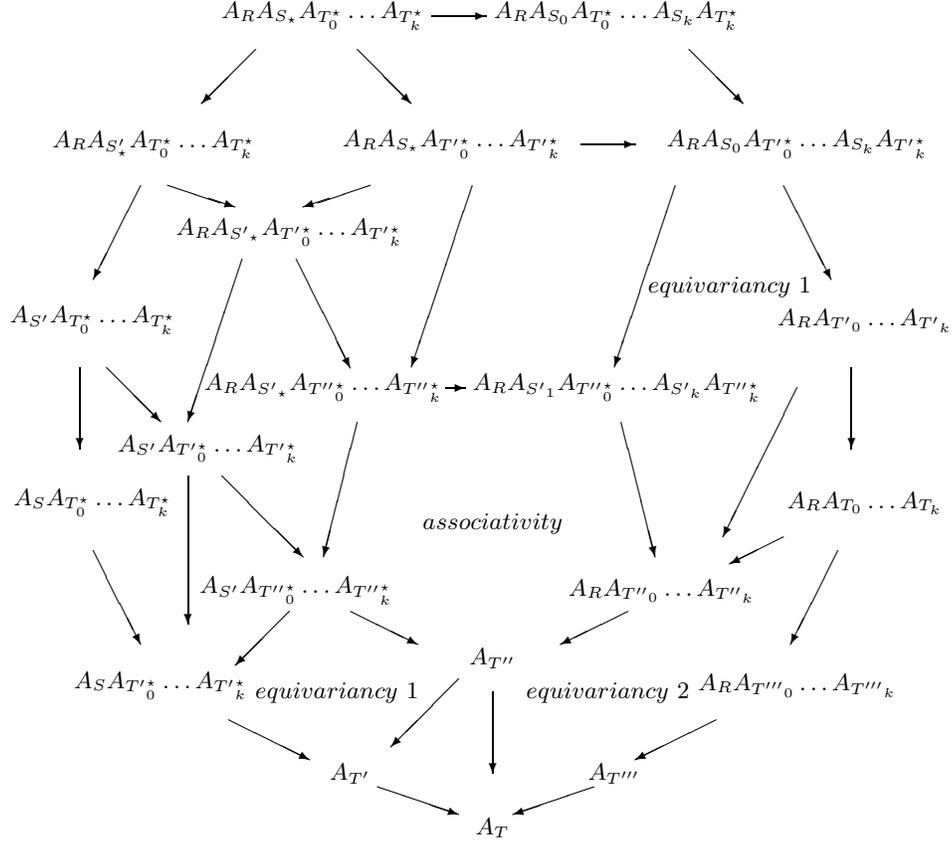
\begin{picture}(300,130)

\put(40,130){\makebox(0,0){\small\mbox{$ A_R
A_{S_{\star}} A_{T_0^{\star}} 
\ldots A_{T_k^{\star}}   
$}}}
\put(30,125){\vector(-1,-1){8}}

\put(56,130){\vector(1,0){8}}

\put(94,125){\vector(1,-1){8}}

\put(84,130){\makebox(0,0){\small\mbox{$ A_R
A_{S_{0}} A_{T_0^{\star}} \ldots
 A_{S_{k}}
A_{T_k^{\star}}   
$ }}}

\put(52,130){\makebox(0,0){\mbox{$ $}}}


\put(15,111){\makebox(0,0){\mbox{\small$A_R
A_{S'_{\star}} A_{T_0^{\star}} 
\ldots A_{T_k^{\star}}   
$}}}

\put(110,111){\makebox(0,0){\mbox{\small$ A_R
A_{S_{0}} A_{{T'}_0^{\star}} \ldots
 A_{S_{k}}
A_{{T'}_k^{\star}}   
$}}}

\put(13,105){\vector(-1,-2){7}}

\put(108,105){\vector(1,-2){7}}


\put(6,85){\makebox(0,0){\mbox{\small$
A_{S'} A_{T_0^{\star}} 
\ldots A_{T_k^{\star}}  
$}}}

\put(120,85){\makebox(0,0){\mbox{\small$A_R
 A_{{T'}_0} \ldots
A_{{T'}_k}   
$}}}

\put(6,58){\makebox(0,0){\mbox{\small$
A_{S} A_{T_0^{\star}} 
\ldots A_{T_k^{\star}}  
$}}}

\put(4,78){\vector(0,-1){12}}

\put(118,78){\vector(0,-1){12}}

\put(120,58){\makebox(0,0){\mbox{\small$A_R
 A_{{T}_0} \ldots
A_{{T}_k}   
$}}}

\put(16,31){\makebox(0,0){\mbox{\small$
A_{S} A_{{T'}_0^{\star}} 
\ldots A_{{T'}_k^{\star}}  
$}}}

\put(6,51){\vector(1,-2){7}}

\put(116,51){\vector(-1,-2){7}}

\put(110,31){\makebox(0,0){\mbox{\small$A_R
 A_{{T'''}_0} \ldots
A_{{T'''}_k}   
$}}}



\put(65,35){\makebox(0,0){\mbox{\small$A_{T''} 
$}}}



\put(44,18){\makebox(0,0){\mbox{\small$A_{T'} 
$}}}

\put(83,18){\makebox(0,0){\mbox{\small$A_{T'''} 
$}}}

\put(65,10){\makebox(0,0){\mbox{\small$A_T 
$}}}


\put(48,16){\vector(3,-1){12}}
\put(80,16){\vector(-3,-1){12}}

\put(26,26){\vector(2,-1){12}}
\put(98,26){\vector(-2,-1){12}}

\put(60,32){\vector(-1,-1){10}}
\put(65,30){\vector(0,-1){12}}

\put(65,55){\makebox(0,0){\small\mbox{$associativity
$}}}

\put(40,75){\makebox(0,0){\small\mbox{$ A_R
A_{{S'}_{\star}} A_{{T''}_0^{\star}} 
\ldots A_{{T''}_k^{\star}}   
$}}}
\put(45,70){\vector(-1,-4){5}}

\put(58,75){\vector(1,0){3}}

\put(84,70){\vector(1,-4){5}}

\put(84,75){\makebox(0,0){\small\mbox{$ A_R
A_{{S'}_{1}} A_{{T''}_0^{\star}} \ldots
 A_{{S'}_{k}}
A_{{T''}_k^{\star}}   
$ }}}


\put(36,45){\makebox(0,0){\mbox{\small$
A_{S'} A_{{T''}_0^{\star}} 
\ldots A_{{T''}_k^{\star}}  
$}}}

\put(44,42){\vector(2,-1){10}}

\put(85,42){\vector(-2,-1){10}}

\put(90,45){\makebox(0,0){\mbox{\small$A_R
 A_{{T''}_0} \ldots
A_{{T''}_k}   
$}}}

\put(23,66){\makebox(0,0){\mbox{\small$
A_{S'} A_{{T'}_0^{\star}} 
\ldots A_{{T'}_k^{\star}}  
$}}}
\put(35,42){\vector(-1,-1){8}}

\put(8,78){\vector(1,-1){8}}

\put(20,62){\vector(0,-1){22}}

\put(25,62){\vector(1,-1){12}}

\put(35,98){\makebox(0,0){\mbox{\small$A_R
A_{{S'}_{\star}} A_{{T'}_0^{\star}} 
\ldots A_{{T'}_k^{\star}}   
$}}}

\put(17,105){\vector(3,-1){10}}
\put(28,94){\vector(-1,-3){8}}
\put(36,94){\vector(1,-2){8}}

\put(42,30){\makebox(0,0){\small\mbox{$equivariancy \ 1
$}}}

\put(82,30){\makebox(0,0){\small\mbox{$equivariancy \ 2
$}}}

\put(100,90){\makebox(0,0){\small\mbox{$equivariancy \ 1
$}}}


\put(59,111){\makebox(0,0){\mbox{\small$A_R
A_{S_{\star}} A_{{T'}_0^{\star}} 
\ldots A_{{T'}_k^{\star}}   
$}}}
\put(62,105){\vector(-1,-3){9}}
\put(45,125){\vector(1,-1){8}}
\put(47,105){\vector(-3,-1){10}}
\put(92,105){\vector(-1,-3){9}}
\put(78,111){\vector(1,0){8}}
\put(110,75){\vector(-1,-2){11}}
\put(108,53){\vector(-2,-1){8}}







\end{picture}}

\noindent In this diagram we omit the symbol $\otimes$ to shorten the notations. Then the central region of the diagram commutes because of associativity of
$A$ with respect to order preserving  maps of $n$-ordinals. Other regions commute
either by one of equivariancy conditions either by naturality either by
functoriality. The commutativity of this diagram means the associativity
of $A$ with respect to composition of maps of $n$-ordinals.    

\Q

\section{The category of quasibijections and configuration spaces.}\label{qconf} 

It is clear that the category $Q_n$ is the union of connected
 components $Q_n(k)$ where $k$ is the cardinality of the $n$-ordinals.  

\begin{theorem}\begin{itemize}\item For a finite $n$ the space  $N(Q_n(k))$    has homotopy type of unordered configuration spaces of $k$-points in $\Re^n ;$ 
\item The  localisation functors
$$l_2:Q_2 \rightarrow \Pi Q_2,$$
induces a weak equivalence of the nerves;
\item The    groupoid  $\Pi  Q_2$ is equivalent to the groupoid of  braids; 
\item The  localisation functors
$$l_{\infty}:Q_{\infty} \rightarrow \Pi Q_{\infty},$$
induces a weak equivalence of the nerves;
\item  the    groupoids   $\Pi  Q_n \ , \ 3\le n\le \infty$ are equivalent to the symmetric groups groupoid. 
\end{itemize}
\end{theorem}

\Proof We give a sketch of the proof. A detailed discussion can be found in \cite{SymBat,Berger}. Consider the configuration space of 
ordered $k$-points in $\Re^n :$
  $$\Conf_k(\Re^n)= \{(x_1,\ldots, x_k)\in (\Re^n)^k \ | \  x_i \ne x_j  \  \mbox{if} \  i\ne j \ \}$$
It admits a so called Fox-Neuwirth stratification.

Let $\stackrel{\scriptscriptstyle o \  \ \ \ \ \ \  }{S^{n-p-1}_{+}}$ denote the  open $(n-p-1)$-hemisphere in $\Re^n$ , $0\le p \le n-1$:
\[  \stackrel{\scriptscriptstyle o \  \ \ \ \ \ \  }{S^{n-p-1}_{+}} = \left\{x\in \Re^n \left \vert \ \begin{array}[c]{ll}
  x_1^2+\ldots+x_n^2 = 1&   \\
   x_{p+1}> 0 \ \mbox{and} \ x_i = 0  \ \mbox{if} \  1\le i \le p &
  \end{array} \right. \right\} \]

Similarly,
 \[  \stackrel{\scriptscriptstyle o \  \ \ \ \ \ \  }{S^{n-p-1}_{-}} = \left\{ x \in \Re^n \left \vert \ \begin{array}[c]{ll}
  x_1^2+\ldots+x_n^2 = 1&   \\
   x_{p+1}< 0 \ \mbox{and} \ x_i = 0  \ \mbox{if} \  1\le i \le p &
  \end{array} \right. \right\} \raisebox{-3mm}{ \ .} \]
Let $u_{ij}: \Conf_k(\Re^n)\rightarrow S^{n-1}$ be the function
$$u_{ij}(x_1,\ldots,x_k) =  \frac{x_j - x_i}{||x_j-x_i||} $$

{ The Fox-Neuwirth cell corresponding to an $n$-ordinal $T$ with $[T]=[k-1]$} is a subspace of $\Conf_{k}(\Re^n)$
\[ FN_T =  \left\{ x\in \Conf_{k}(\Re^n)  \left |
 \ \begin{array}[c]{ll}
  \ u_{ij}(x)  \in \stackrel{\scriptscriptstyle o \  \ \ \ \ \ \  }{S^{n-p-1}_{+}}&  \ \mbox{\rm if} \ i<_p j  \ \mbox{\rm in} \ T  \\
  &  \\
  \ u_{ij}(x)  \in \stackrel{\scriptscriptstyle o \  \ \ \ \ \ \  }{S^{n-p-1}_{-}}&  \ \mbox{\rm if} \ j<_p i  \ \mbox{\rm in} \ T
  \end{array}\right.  \right\} \raisebox{-7mm}{ \ .}\]

Each Fox-Neuwirth cell is an open convex subspace of $(\Re^n)^k .$ We also have 
 $$\Conf_k(\Re^n)= \bigcup_{[T]= [k-1] \ , \ \pi \in S_k} \pi FN_T.$$ 
Here $\pi FN_T$ means a space obtained from $FN_T$ by renumbering  points according to the permutation $\pi .$ 

Let $J_n(k)$ be the Milgram poset of all possible $n$-ordinal structures on the set $\{0,\ldots,k-1\} $ \cite{SymBat}. 
The group $S_k$ acts  on $J_n(k)$ and the quotient $J_n(k)/S_k$ is isomorphic to $Q_n(k).$ 

One can think of an element from $J_n(k)$ as a pair $(T,\pi)$ where $T$ is an $n$-ordinal and $\pi$ is a permutation from $S_k $ and  
 $(T,\pi) > (S,\xi)$ in $ J_n(k)$ when there exists a quasibijection $\sigma:T\rightarrow S$ and $\xi\cdot\pi = \sigma .$ 

We also can associate a convex subspace of the configuration space $FN(T,\pi) = \pi FN_T$ with every element of $J_n(k) .$ Moreover,
if $(T,\pi) > (S,\xi)$ then $FN(S,\xi)$ is on the boundary of the closure of  $FN(T,\pi).$ Let us define 
$$\overline{FN}(T,\pi)= \bigcup_{(S,\xi)\le (T,\pi)} FN(S,\xi).$$
The spaces $\overline{FN}(T,\pi)$ are contractible and, moreover, we have a functor $$\overline{FN}: J^{op}_n(k)\rightarrow Top .$$
We then have the following zig-zag of weak equivalences
$$N(J_n^{op}(k))\leftarrow hocolim~ \overline{FN} \rightarrow colim~ \overline{FN} \simeq  \Conf_k(\Re^n) .$$
The first statement of the theorem follows then from the quotient of the zig-zag above by the action of the symmetric group.  
 The second and the third statements are the consequences of the fact that
the space  $\Conf_k(\Re^2)$ is the $K(Br_k,1)$-space.
The fifth   statement follows from the fact that  the fundamental group of $\Conf_k(\Re^n)$ is trivial for $  n>3 .$ Finally the fourth statement can be obtained using the formula $Q_{\infty} = colim_n Q_n .$  
 
\
    
\Q

\



            We shall now, in Lemmas \ref{5.1} and \ref{spliting}, make the equivalence
between $\Pi Q_2$ and $\Br$ more explicit. These results will then be used in
section \ref{vs} to relate different operadic notions.

The total order functor $[-]:Q_2\rightarrow \Sm$ induces by the universal property a functor $s_2:\Pi Q_2 \rightarrow \Sm .$
 
Let $p:\Br\rightarrow \Sm$ be the canonical functor.
The map $p$ admits a section $q ,$ which is not  a homomorphism. For $\sigma\in S_n$ we construct a braid
$q(\sigma)$ which for $i<j$ such that $\sigma(i)> \sigma(j)$ has a strand from $i$ to $\sigma(i)$ which goes over the strand from $j$ to $\sigma(j) $ and there is no crossing if $\sigma$ preserves the order of $i$ and $j .$
  
\begin{lem}\label{5.1} \begin{itemize}
\item The composite
$$Q_2\stackrel{[-]}{\longrightarrow} \Sm\stackrel{q}{\longrightarrow} \Br$$
is a functor; 
\item The  functor induced by the universal property of $\Pi Q_2$
$$b:\Pi Q_2 \rightarrow \Br$$
is an equivalence of groupoids.
\end{itemize}
\end{lem} 

\Proof To prove that $q[-]$ is a functor we have to prove that it preserves composition. We observe that in a composite of quasibijections of $2$-ordinals $T\stackrel{\sigma}{\rightarrow}S\stackrel{\xi}{\rightarrow} R$
if $\sigma$ reverses the total order of two elements $i,j\in T$  then $\xi$ can not reverse the order of $\sigma(i)$ and $\sigma(j).$ So, the resulting overcrossings in the composite
$q[\sigma]q[\xi]$ are the same as in $q[\sigma\cdot\xi] .$  

To prove the second claim it is sufficient to check that the induced morphism of groups
$$b:\Pi Q_2(S[n-1],S[n-1])\rightarrow \Br_n$$
is an isomorphism. 

It is obviously an epimorphism. So we have to prove that it is also a monomorphism.

For this it will be enough to prove that if a zig-zag
$$z:S[n-1]\leftarrow T[n-1]\rightarrow S[n-1] \leftarrow \ldots \rightarrow S[n-1],$$
where each arrow is given by a permutation of two consecutive elements or an identity permutation, is such that the corresponding braid $b(z)$ is trivial then  $z$ is trivial in $\Pi Q_2 .$ 

This can be done if we prove that the morphisms in $\Pi Q_2(S[n-1],S[n-1])$
$$\bar{\sigma}_i: S[n-1]\stackrel{1}{\longleftarrow} T[n-1]\stackrel{\sigma_i}{\longrightarrow} S[n-1],$$  
where the left arrow is given by an identity and the right arrow is given by permutation $\sigma_i$ which change the order of $i$ and $i+1,$ satisfy the classical Artin braid relations. Then we can   prove   triviality of $z$ using the same rewriting process as for $b(z) .$

Let $j>i+1$ and choose $m,l$ such that 
$[m-1]\oplus[l-1] = [n-1]$ and $i\in [m-1] , j=m+1.$ 
The following commutative diagram in $Q_2$ proves that $\bar{\sigma}_i\bar{\sigma}_j= \bar{\sigma}_j\bar{\sigma}_i :$ 

 {\unitlength=1mm

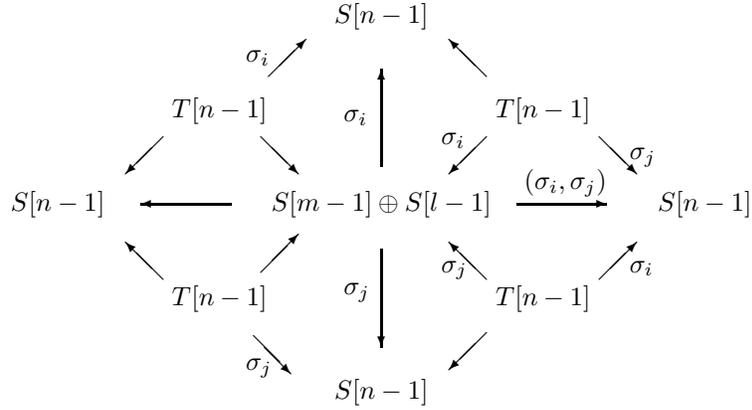
\begin{picture}(200,60)(-10,-5)

\put(0,25){\makebox(0,0){\mbox{$S[n-1]$}}}

\put(23,25){\vector(-1,0){12}}
\put(72,15){\vector(1,1){5}}
\put(72,34){\vector(1,-1){5}}

\put(14,34){\vector(-1,-1){5}}

\put(14,15){\vector(-1,1){5}}
\put(26,7.5){\vector(1,-1){5}}

\put(57,8){\vector(-1,-1){5}}
\put(57,15){\vector(-1,1){5}}
\put(27,16){\vector(1,1){5}}
\put(28,42){\vector(1,1){5}}
\put(57,42){\vector(-1,1){5}}
\put(27,34){\vector(1,-1){5}}
\put(57,34){\vector(-1,-1){5}}

\put(43,30){\vector(0,1){13}}
\put(43,19){\vector(0,-1){13}}

\put(43,25){\makebox(0,0){\mbox{$S[m-1]\oplus S[l-1]$}}}
\put(43,50){\makebox(0,0){\mbox{$S[n-1]$}}}
\put(43,0){\makebox(0,0){\mbox{$S[n-1]$}}}

\put(21.5,12.5){\makebox(0,0){\mbox{$T[n-1]$}}}
\put(64.5,12.5){\makebox(0,0){\mbox{$T[n-1]$}}}

\put(21.5,37.5){\makebox(0,0){\mbox{$T[n-1]$}}}
\put(64.5,37.5){\makebox(0,0){\mbox{$T[n-1]$}}}

\put(61,25){\vector(1,0){12}}

\put(86,25){\makebox(0,0){\mbox{$S[n-1]$}}}
\put(25,44){\shortstack{\mbox{$\sigma_i $}}}
\put(38,36){\shortstack{\mbox{$\sigma_i $}}}
\put(51,33){\shortstack{\mbox{$\sigma_i $}}}
\put(76,16){\shortstack{\mbox{$\sigma_i $}}}

\put(25,3){\shortstack{\mbox{$\sigma_j $}}}
\put(38,13){\shortstack{\mbox{$\sigma_j $}}}
\put(51,16){\shortstack{\mbox{$\sigma_j $}}}
\put(76,31){\shortstack{\mbox{$\sigma_j $}}}

\put(62,27){\shortstack{\mbox{$(\sigma_i,\sigma_j) $}}}
\end{picture}}

In this diagram all unnamed morphisms are identities on the underlying sets. The morphism $(\sigma_i,\sigma_j)$ acts as $\sigma_i$ on $[m-1]$ and as $\sigma_0$ on $[l-1] .$ 
    
For the proof of Yang-Baxter relations $\bar{\sigma}_i\bar{\sigma}_{i+1}\bar{\sigma}_i = \bar{\sigma}_{i+1}\bar{\sigma}_i\bar{\sigma}_{i+1} $ 
we should consider the following commutative diagram in 
$Q_2 $ which expresses the morphism $\bar{\sigma}_{i+1}\bar{\sigma}_i\bar{\sigma}_{i+1}$

 {\unitlength=1mm

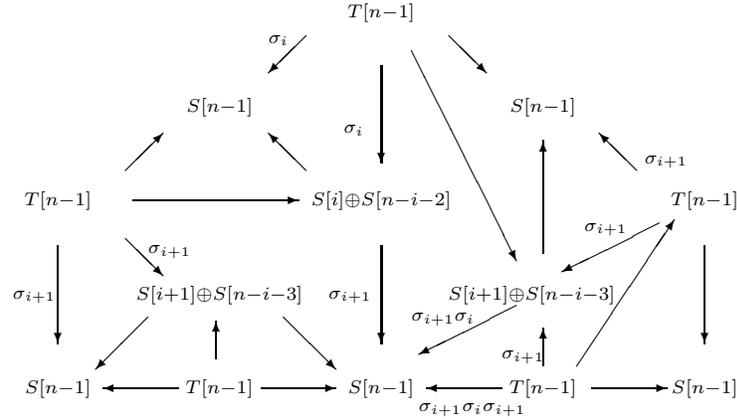
\begin{picture}(200,60)(-10,-5)

\put(0,25){\makebox(0,0){\mbox{$\scriptstyle T[n-1]$}}}

\put(10,25){\vector(1,0){22}}
\put(77,29){\vector(-1,1){5}}

\put(9,29){\vector(1,1){5}}

\put(9,20){\vector(1,-1){5}}
\put(30,9.5){\vector(1,-1){7}}

\put(12,9.5){\vector(-1,-1){7}}
\put(33,47){\vector(-1,-1){5}}
\put(52,47){\vector(1,-1){5}}
\put(33,29){\vector(-1,1){5}}

\put(64.5,3){\vector(0,1){5}}
\put(69,3){\vector(2,3){13}}
\put(80,22){\vector(-2,-1){13}}
\put(61,11){\vector(-2,-1){13}}
\put(64.5,18){\vector(0,1){15}}
\put(47,45){\vector(1,-2){14}}

\put(43,43){\vector(0,-1){13}}
\put(43,19){\vector(0,-1){13}}
\put(86,19){\vector(0,-1){13}}
\put(0,19){\vector(0,-1){13}}

\put(43,25){\makebox(0,0){\mbox{$\scriptstyle S[i]\oplus S[n-i-2]$}}}
\put(43,50){\makebox(0,0){\mbox{$\scriptstyle T[n-1]$}}}
\put(43,0){\makebox(0,0){\mbox{$\scriptstyle S[n-1]$}}}
\put(21.5,0){\makebox(0,0){\mbox{$\scriptstyle T[n-1]$}}}
\put(64.5,0){\makebox(0,0){\mbox{$\scriptstyle T[n-1]$}}}
\put(0,0){\makebox(0,0){\mbox{$\scriptstyle S[n-1]$}}}
\put(86,0){\makebox(0,0){\mbox{$\scriptstyle S[n-1]$}}}

\put(21.5,12.5){\makebox(0,0){\mbox{$\scriptstyle S[i+1]\oplus S[n-i-3]$}}}
\put(63,12.5){\makebox(0,0){\mbox{$\scriptstyle S[i+1]\oplus S[n-i-3]$}}}

\put(21.5,37.5){\makebox(0,0){\mbox{$\scriptstyle S[n-1]$}}}
\put(64.5,37.5){\makebox(0,0){\mbox{$\scriptstyle S[n-1]$}}}
\put(21,4){\vector(0,1){5}}
\put(71,0){\vector(1,0){10}}
\put(27,0){\vector(1,0){10}}
\put(16,0){\vector(-1,0){10}}
\put(59,0){\vector(-1,0){10}}

\put(86,25){\makebox(0,0){\mbox{$\scriptstyle T[n-1]$}}}

\put(59,4){\shortstack{\mbox{$\scriptstyle \sigma_{i+1} $}}}
\put(70,21){\shortstack{\mbox{$\scriptstyle \sigma_{i+1} $}}}
\put(78,30){\shortstack{\mbox{$\scriptstyle \sigma_{i+1} $}}}
\put(36,12){\shortstack{\mbox{$\scriptstyle \sigma_{i+1} $}}}
\put(38,34){\shortstack{\mbox{$\scriptstyle \sigma_{i} $}}}
\put(-6,12){\shortstack{\mbox{$\scriptstyle \sigma_{i+1} $}}}
\put(12,18){\shortstack{\mbox{$\scriptstyle \sigma_{i+1} $}}}
\put(28,46){\shortstack{\mbox{$\scriptstyle \sigma_{i} $}}}
\put(47,9){\shortstack{\mbox{$\scriptstyle \sigma_{i+1}\sigma_{i} $}}}
\put(48,-3){\shortstack{\mbox{$\scriptstyle \sigma_{i+1}\sigma_{i}\sigma_{i+1} $}}}

\end{picture}}

An analogous diagram (the mirror image of the above diagram) can be written for $\bar{\sigma}_{i}\bar{\sigma}_{i+1}\bar{\sigma}_{i}.$
The relation follows from it immediately.

\Q 

So, we have a commutative diagram of categories and functors

 {\unitlength=1mm
\begin{picture}(200,30)(-10,5)
\put(20,25){\makebox(0,0){\mbox{$Q_2$}}}
\put(24,25){\vector(1,0){13}}
\put(20,22){\vector(2,-1){16}}
\put(64,22){\vector(-2,-1){16}}
\put(43,25){\makebox(0,0){\mbox{$\Pi Q_2$}}}
\put(43,22){\vector(0,-1){7}}
\put(43,10){\makebox(0,0){\mbox{$\Sm$}}}
\put(49,26){\vector(1,0){13}}
\put(62,24){\vector(-1,0){13}}
\put(54,21){\shortstack{\mbox{$c $}}}
\put(19,17){\shortstack{\mbox{$[-] $}}}
\put(38,17){\shortstack{\mbox{$s_2 $}}}
\put(58,17){\shortstack{\mbox{$p $}}}
\put(66,25){\makebox(0,0){\mbox{$\Br$}}}
\put(54,27){\shortstack{\mbox{$b $}}}
\end{picture}}

\noindent where $c$ is an adjoint equivalence to $b .$ Notice that  all functors in this diagram are strict monoidal functors. 


 \begin{lem}\label{spliting} Let $$z: S\stackrel{\sigma}{\longleftarrow} T \stackrel{\eta}{\longrightarrow} R$$ be a zig-zag of quasibijections of $n$-ordinals such that $$s_2(z) = \tau_1\oplus \ldots\oplus \tau_k .$$ Then there exist braids $ b_i \ , \ 1\le i\le k$ such that $p(b_i) = \tau_i \ , \ 1\le i\le k$ and 
$$b(z) = b_1\oplus\ldots\oplus b_k .$$

 \end{lem}   
      
\Proof We will prove that there exist  quasibijections $\sigma_{i}:T_i\rightarrow S_i = S[n_i]  \ , \ \eta_{i}:T_i\rightarrow R_i= S[n_i] \ 1\le i \le k , $ two quasibijections $\xi:\oplus_i S_i \rightarrow S \ , \ \zeta:\oplus_i R_i \rightarrow R$ and a quasibijection 
$\kappa: \oplus_i T_i\rightarrow T ,$  such that the following diagram commutes

{\unitlength=1mm

\begin{picture}(200,25)(-10,5)
\put(66,25){\makebox(0,0){\mbox{$\oplus_i R_i$}}}
\put(20,25){\makebox(0,0){\mbox{$\oplus_i S_i$}}}
\put(20,10){\makebox(0,0){\mbox{$S$}}}

\put(37,10){\vector(-1,0){13}}
\put(20,20){\vector(0,-1){6}}
\put(37,25){\vector(-1,-0){12}}
\put(66,20){\vector(0,-1){6}}
\put(48,25){\vector(1,0){12}}

\put(43,25){\makebox(0,0){\mbox{$\oplus_i T_i$}}}
\put(43,20){\vector(0,-1){6}}
\put(40,17){\shortstack{\mbox{$\kappa $}}}
\put(43,10){\makebox(0,0){\mbox{$T$}}}
\put(17,17){\shortstack{\mbox{$\xi $}}}\put(63,17){\shortstack{\mbox{$\zeta $}}}
\put(49,10){\vector(1,0){13}}

\put(66,10){\makebox(0,0){\mbox{$R$}}}

\put(30,11){\shortstack{\mbox{$\sigma $}}}
\put(54,11){\shortstack{\mbox{$\eta $}}}
\put(28,26){\shortstack{\mbox{$\oplus_i\sigma_i $}}}
\put(51,26){\shortstack{\mbox{$\oplus_i \eta_i $}}}
\end{picture}}  
\noindent and $b(\xi) = b(\zeta) = \Gamma_B(\pi;1,\ldots,1)$ for a braid $\pi$ on $k$ strands. 
Then the result will follow from an elementary observation that the braid $$b(S)\stackrel{b(\xi)^{-1}}{\longrightarrow} \oplus_i b(S_i)   
\stackrel{\oplus_i b(\sigma_i)^{-1}}{\longrightarrow}\oplus_i b(T_i)\stackrel{\oplus_i b(\eta_i^{})}{\longrightarrow}\oplus_i b(R_i)\stackrel{b(\xi)}{\longrightarrow} b(R)$$
is equal to
$$ \oplus_i b(S_i)   
\stackrel{\oplus_i b(\sigma_i^{-1})}{\longrightarrow}\oplus_i b(T_i)\stackrel{\oplus_i b(\eta_i^{})}{\longrightarrow}\oplus_i b(R_i) .$$

It is enough to proof  the lemma for $k=2 .$ The rest will follow by induction. Also without loss of generality we can assume that 
$S = S[n]$ and $T=T[n] .$ Now, $p(S)$ is the ordinal sum $[l]\oplus [m], n= m+1+1 $ and the image of the restriction of the map $\sigma^{-1}\eta$ on $\{0,\ldots ,l\}$ is $\{0,\ldots ,l\}$ and the image of the restriction on $\{l+1,\ldots ,m+l+1\}$
is $\{l+1,\ldots ,m+l+1\}.$ 

We put $S_1 = S[l], T_1 = T[l]$ and $S_2 = S[m], T_1 = T[m] .$ 
We have to construct quasibijections $$\sigma_i,\eta_i:T_i\rightarrow S_i \ i=1,2 ,$$
and also quasibijections 
$$\xi,\zeta: S_1\oplus S_2 \rightarrow S \ , \kappa: T_1 \oplus T_2 \rightarrow T ,$$
which make  the diagram \newpage
 $$  S_1\oplus S_2\stackrel{\sigma_1 \oplus \sigma_2}{\longleftarrow} T_1\oplus T_2 \stackrel{\eta_1\oplus \eta_2}{\longrightarrow} S_1\oplus S_2$$   
\begin{equation}\label{ds}\downarrow \hspace{20mm} \downarrow \hspace{20mm} \downarrow\end{equation}
$$S\hspace{7mm}\stackrel{\sigma}{\longleftarrow}\hspace{5mm} T \hspace{5mm}\stackrel{\eta}{\longrightarrow}\hspace{7mm} S$$
 commutative.

The quasibijection $\kappa$ is simply the identity. 
Let us describe $\sigma_1 .$ Let $\sigma([l])$ be the image of the set $\{0,\ldots,l\}$ in the ordinal $[n] .$ This image gets an induced order from $[n] $ which makes it isomorphic to $[l] .$
Let $\phi_1:\sigma([l])\rightarrow [l]$ be this unique isomorphism. We define $\sigma_1$ as the composite
$$[l]\rightarrow \sigma([l])\stackrel{\phi_1}{\rightarrow} [l] .$$ 
Similarly, we define $\sigma_2$ as the composite 
$$[m]\rightarrow \sigma([m])\stackrel{\phi_2}\rightarrow [m] , $$
where $\sigma([m])$ is the image of $\{l+1,\ldots,m+l+1\}$ 
 and we give analogous definitions for $\eta_{1}$ and $\eta_{2} .$ 

Finally, we   define $\xi$ by the formula
$$\xi(x) =  \left\{
\begin{array}{rl}
\phi_1^{-1}(x) & \mbox{if}  \ \ \ \   x\in \{0,\ldots,l\} \\
\phi_2^{-1}(x) & \mbox{if} \ \ \ \ \ x\in \{l+1,\ldots,m+l+1\}
\end{array} \right. $$    
We use a similar argument to define $\zeta.$ The commutativity of the diagram (\ref{ds}) follows from the definition.  

 \Q

\section{Quasisymmetric $n$-operads vs symmetric and braided operads.}\label{vs}

\begin{theorem}\label{qvsb} The category of quasisymmetric $2$-operads and the category of braided operads are equivalent. 
 \end{theorem}
 
\Proof 
We first prove that the category of quasisymmetric $2$-operads is equivalent to the category whose objects are   {\it mixed $2$-operads} in the sense of the definition below and whose morphisms are multiplications and units preserving morphisms of the underlying braided collections.

 \begin{defin}\label{defbrop2} A mixed  $2$-operad in $V$ is a right braided collection $A$ equipped with the following additional structure:

- a morphism  $e:I\rightarrow A_0$

- for every order preserving  map $\sigma:[n]\rightarrow [k]$ in $\Omega^s$   a morphism
: 
$$\mu_{\sigma}: A_{k}\otimes(A_{n_0}\otimes ... \otimes
A_{n_k})\longrightarrow A_{n},
$$
where $[n_i] = \sigma^{-1}(i).$ 

They must satisfy the   identities (1-3) from the definition of symmetric operad  
and  the following 
  two equivariance conditions:

\begin{enumerate}
\item

For any two quasibijections of $2$-ordinals   $\pi, \rho$ and two order preserving maps $\sigma,\sigma'\in \Omega^s$ such that 
the following  diagram commutes in $\Omega^s $

{\unitlength=0.7mm

\begin{picture}(40,28)(-45,3)

\put(13,25){\makebox(0,0){\mbox{$ [T']$}}}
\put(13,21){\vector(0,-1){10}}
\put(5,15){\shortstack{\mbox{$[\pi] $}}}
\put(43,15){\shortstack{\mbox{$[\rho] $}}}

\put(22,25){\vector(1,0){10}}

\put(26,26){\shortstack{\mbox{$\sigma' $}}}

\put(41,25){\makebox(0,0){\mbox{$ [S']$}}}
\put(41,21){\vector(0,-1){10}}

\put(13,7){\makebox(0,0){\mbox{$ [T]$}}}

\put(23,7){\vector(1,0){10}}

\put(27,8){\shortstack{\mbox{$\sigma $}}}

\put(41,7){\makebox(0,0){\mbox{$[S]$}}}

\end{picture}}

the following induced diagram commutes:

{\unitlength=1mm

\begin{picture}(200,33)(-30,0)

\put(9,25){\makebox(0,0){\mbox{$A_{k'}\otimes(A_{n'_{\rho(0)}}
\otimes ... \otimes
A_{n'_{\rho(k)}}) $}}}
\put(10,10){\vector(0,1){10}}
\put(-9,15){\shortstack{\mbox{$\scriptstyle A(b(\rho))\otimes\tau(\rho) $}}}

\put(32,25){\vector(1,0){12}}

\put(35,27){\shortstack{\mbox{$\mu_{\sigma'} $}}}

\put(50,25){\makebox(0,0){\mbox{$ A_{n'}$}}}
\put(50,10){\vector(0,1){10}}

\put(10,5){\makebox(0,0){\mbox{$A_{k}\otimes(A_{n_0}
\otimes ... \otimes
A_{n_k})$}}}

\put(31,5){\vector(1,0){13}}

\put(35,7){\shortstack{\mbox{$\mu_{\sigma} $}}}

\put(50,5){\makebox(0,0){\mbox{$A_{n}$}}}

\put(56,15){\makebox(0,0){\mbox{$\scriptstyle A(b(\pi)) $}}}

\put(60,4){\shortstack{\mbox{$, $}}}

\end{picture}}
\noindent where $\tau(\rho)$ is the symmetry in $V$  which corresponds to the permutation $[\rho] .$ 

\item For any two quasibijections $\sigma, \sigma'$ and two order preserving maps $\eta,\eta'\in \Omega^s$ such that 
the following  diagram commutes in $\Omega^s $

{\unitlength=0.9mm
\begin{picture}(40,29)(-29,3)

\put(13,25){\makebox(0,0){\mbox{$ [T'']$}}}
\put(13,21){\vector(0,-1){10}}
\put(6,15){\shortstack{\mbox{$[\sigma] $}}}
\put(42,15){\shortstack{\mbox{$\eta' $}}}
\put(22,25){\vector(1,0){10}}
\put(24,27){\shortstack{\mbox{$[\sigma'] $}}}
\put(41,25){\makebox(0,0){\mbox{$ [T']$}}}
\put(41,21){\vector(0,-1){10}}
\put(13,7){\makebox(0,0){\mbox{$ [T]$}}}
\put(23,7){\vector(1,0){10}}
\put(27,8.5){\shortstack{\mbox{$\eta $}}}
\put(41,7){\makebox(0,0){\mbox{$[S]$}}}
\end{picture}}

\noindent  the following diagram commutes 

{\unitlength=1mm

\begin{picture}(200,53)(-30,0)

\put(10,25){\makebox(0,0){\mbox{$A_{k}\otimes (A_{n''_0}
\otimes ... \otimes
A_{n''_k}) $}}}
\put(10,8){\vector(0,1){12}}
\put(-23,13){\shortstack{\mbox{$\scriptstyle 1\otimes A(b(\sigma_0))\otimes\ldots\otimes A(b(\sigma_k)) $}}}

\put(-23,34){\shortstack{\mbox{$\scriptstyle 1\otimes A(b(\sigma'_0))\otimes\ldots\otimes A(b(\sigma'_k)) $}}}

\put(46,34){\shortstack{\mbox{$\scriptstyle  A(b(\sigma')) $}}}
\put(46,13){\shortstack{\mbox{$\scriptstyle  A(b(\sigma)) $}}}
\put(35,47){\makebox(0,0){\mbox{$ \mu_{\eta'}$}}}
\put(34,7){\makebox(0,0){\mbox{$ \mu_{\eta}$}}}

\put(45,25){\makebox(0,0){\mbox{$ A_{n''}$}}}
\put(45,8){\vector(0,1){12}}

\put(75,25){\makebox(0,0){\mbox{$ $}}}

\put(57,21){\shortstack{\mbox{$ $}}}

\put(10,5){\makebox(0,0){\mbox{$A_{k}\otimes (A_{n_0}
\otimes ... \otimes
A_{n_k})$}}}

\put(29,5){\vector(1,0){10}}

\put(27,6){\shortstack{\mbox{$ $}}}

\put(45,5){\makebox(0,0){\mbox{$A_{n}$}}}

\put(75,5){\makebox(0,0){\mbox{$ $}}}


\put(60,6){\shortstack{\mbox{$ $}}}


\put(10,45){\makebox(0,0){\mbox{$A_{k}\otimes(A_{n'_0}
\otimes ... \otimes
A_{n'_k}) $}}}

\put(45,45){\makebox(0,0){\mbox{$ A_{n'}$}}}

\put(29,45){\vector(1,0){10}}

\put(10,41){\vector(0,-1){12}}
\put(45,41){\vector(0,-1){12}}

\end{picture}}

\end{enumerate}
\end{defin}  
 For a quasisymmetric $2$-operad $A$ we define a mixed $2$-operad $B$ by pulling back along the equivalence $c:\Br\rightarrow \Pi Q_2 .$ And vice versa, we produce a quasisymmetric $2$-operad from a mixed $2$-operad by pulling back along $b:\Pi Q_2 \rightarrow Br .$ It is not hard to check that this indeed gives the necessary equivalence of the corresponding operadic categories. 

Now, we will prove that the category of mixed $2$-operads is equivalent to the category of braided operads.
Let $A$ be an operad in the sense of \ref{defbrop2}. We have to check that $A$ also satisfies the Fiedorowicz equivariance conditions. Let us start from the second condition.

For each $\rho_i$ let us choose a   zigzag of  $2L$ morphisms in $Q_2 ,$ such that 
$$\rho_i = b(T_i\stackrel{\tau_{1}}{\longleftarrow} R_{i}^1\stackrel{\tau_{2}}{\longrightarrow} R_{i}^2\leftarrow \ldots\longleftarrow  R_{i}^{2L}
\stackrel{\tau_{2k}}{\longrightarrow} S_i ).$$

Obviously, such a zig-zag exists and $L$ can be chosen independently on $i .$ Then the following square commutes for each odd  $j $ :

{\unitlength=0.9mm
\begin{picture}(40,30)(-29,2)

\put(13,25){\makebox(0,0){\mbox{$ [n]$}}}
\put(13,21){\vector(0,-1){10}}
\put(0,15){\shortstack{\mbox{$\scriptstyle [\oplus_i \tau_i^{j}] $}}}
\put(52,15){\shortstack{\mbox{$\sigma $}}}
\put(18,25){\vector(1,0){27}}
\put(24,27){\shortstack{\mbox{$\scriptstyle [\oplus_i \tau_i^{j+1}] $}}}
\put(50,25){\makebox(0,0){\mbox{$ [n]$}}}
\put(50,21){\vector(0,-1){10}}
\put(13,7){\makebox(0,0){\mbox{$ [n]$}}}
\put(18,7){\vector(1,0){27}}
\put(31,8.5){\shortstack{\mbox{$\sigma $}}}
\put(50,7){\makebox(0,0){\mbox{$[k]$}}}
\end{picture}}

\noindent Hence, the application of the second equivariance condition of definition \ref{defbrop2}
$L$ times gives the second Fiedorowicz equivariance condition.

For the first equivariance condition we do an analogous construction by  choosing a presentation of the braid $\rho$ as an image of a zigzag. 

Let $A$ be an operad in the sense of \ref{defbrop}. We construct 
an operad $B$ in the sense of \ref{defbrop2} as follows. As a braided collection $B$ coincides with $A .$ Its multiplication is the same as in $A$ also. The only nontrivial statement to check is that $B$ satisfies the equivariance conditions from  Definition \ref{defbrop2}.    
To prove the second condition we use  Lemma \ref{spliting}.
 
  It is obvious also that the first equivariance condition is satisfied in the following special case. Let $\sigma':T\rightarrow S'$ be an order preserving map and let $\rho:S'\rightarrow S$ be a quasibijection. Apply Lemma \ref{factorisation} to produce a quasibijection $\pi(\rho,\sigma'):T'\rightarrow T$ and order preserving map $\sigma(\rho,\sigma'):T\rightarrow S$ such that $\sigma'\cdot\rho = \pi(\rho,\sigma')\cdot\sigma(\rho,\sigma').$
Then $b(\pi(\rho,\sigma')) = \Gamma_B(b(\rho);1,\ldots,1)$ and we can apply the first equivariance Fiedorowicz condition. 

  Then the first equivariance condition is satisfied in general because of the second equivariance condition of the Definition \ref{defbrop2} applied to the commutative diagram

{\unitlength=0.9mm

\begin{picture}(40,28)(-29,4)

\put(13,25){\makebox(0,0){\mbox{$ [T']$}}}
\put(13,21){\vector(0,-1){10}}
\put(5,15){\shortstack{\mbox{$[\pi] $}}}
\put(47,15){\shortstack{\mbox{$[\sigma] $}}}
\put(18,25){\vector(1,0){22}}
\put(19,27){\shortstack{\mbox{$[\pi(\rho,\sigma'))]$}}}
\put(45,25){\makebox(0,0){\mbox{$ [T]$}}}
\put(45,21){\vector(0,-1){10}}
\put(13,7){\makebox(0,0){\mbox{$ [T]$}}}
\put(18,7){\vector(1,0){22}}
\put(27,9){\shortstack{\mbox{$[\sigma] $}}}
\put(45,7){\makebox(0,0){\mbox{$[S]$}}}

\end{picture}}

\Q

\begin{theorem}\label{qvss} The category of $Q_n$-operads $3\le n\le \infty$ and the category of symmetric
operads are equivalent. \end{theorem}

\Proof The proof is a repetition of the above proof with a simplification that $s_n:\Pi Q_n \rightarrow \Sm$ for $ 3\le n\le \infty$ is an equivalence.

\Q

\section{Locally constant $n$-operads.}\label{lcoperad}

The quasisymmetric $n$-operads are defined in any symmetric monoidal category $V .$ But according to Theorems \ref{qvsb} and \ref{qvss}  they are different from symmetric operads only when $n= 1,2. $ 
As we have seen before the main reason why quasisymmetric operads collapse to symmetric operads for $n>2$ is that configuration space $Conf_k(\Re^n)$ is simply connected and so  localising with respect to quasibijections can only produce a groupoid equivalent to $\Sm .$ The correct procedure, therefore, should be to take the weak $\omega$-groupoid $\Pi_{\infty} Q_n$ and  consider presheaves on it with values in $V$  as the category of collections. There are, however, considerable technical difficulties with this approach. 

Fortunately, the results of Cisinski \cite{cis} show a way around this problem by considering as the category of collections the category of {\it locally constant} functors from $Q_n^{op}$ to $V .$ Pursuing this idea we give the following definition.
\begin{defin}
 Let $V$ be a symmetric monoidal  category and $\W$ (weak equivalences) be a  subclass of its morphisms. A locally constant $n$-operad in $(V,\W)$ is an $n$-operad $A$ in $V$ such that for every quasibijection $\sigma:T\rightarrow S$ the morphism $A(\sigma): A_S\rightarrow A_T$ is a weak equivalence.
\end{defin}

\Remark  We have chosen the name locally constant $n$-operads (which some people prefer to call homotopically locally constant $n$-operads) for two reasons. First, we would like our terminology to agree with the terminology of \cite{cis}. But more important  reason is about philosophy.  The notion of locally constant $n$-operad (and locally constant functor) depends only on the class of weak equivalences but not on the choice of homotopy theory in $V.$ For example,
if $V$ is a symmetric monoidal category and $Iso$ is the class of all isomorphisms, a locally constant $n$-operad in $(V,Iso)$ is the same as a quasisymmetric $n$-operad in $V.$ So, the word `homotopical' is a little bit misleading. Compare this situation with the theory of homotopy limits developed in \cite{Hirch}. We believe that a `true' reason for this phenomenon is that homotopy limit and locally     constant functors are higher categorical rather then homotopical notions. But the homotopy theory is helpful in computations.     As far as we know a similar argument is behind Cisinski's choice of terminology.

An example of an interesting locally constant $n$-operad in the model category of topological spaces, which is not a quasisymmetric $n$-operad is the Getzler-Jones
$n$-operad $GJ^n$ constructed in \cite{SymBat} for all $n<\infty .$ One can also construct an $\infty$-version  $GJ^{\infty}$ by the formula $GJ^{\infty}_T = GJ^n_{\overline{T}},$ where $<_{-n}$ is the minimal  nonempty relation in the $\infty$-ordinal $T,$ the $n$-ordinal $\overline{T}$ has the same underlying set as $T$ and the relation  $<_{n-p-1}$    in $\overline{T}$ coincides with the relation  $<_{-p}$ in $T .$      

 Let $V$ be a symmetric monoidal category equipped with a class of weak equivalences $\W .$   We introduce the following notations:
\begin{itemize}
\item $SO$ is the category of symmetric operads in $V$;
\item $BO$ is the category of braided operads in $V$;
\item $O_n$ is the category of $n$-operads in $V$;
\item $QO_n$ is the full subcategory of $O_n$ of quasisymmetric $n$-operads in $V$;
\item $LCO_n$ is the full subcategory of $O_n$ of  locally constant $n$-operads in $(V,\W)$.
\end{itemize}

\begin{defin} A morphism of operads (in any of the categories above) is a weak equivalence if it is a termvise weak equivalence of the collections. The homotopy category of operads is the category of operads localised with respect to the class of weak equivalences.\end{defin}

Let us describe the relations between the different categories of operads we deal with in this paper. We have  already done it for the case $\W=Iso$ in Section \ref{vs}.

Let us fix a base symmetric monoidal model category $V$ and let $\W$ be its class of weak equivalences
in the model category theoretic sense.
Moreover, we will assume that $V$ satisfies the conditions from Section 5 of \cite{SymBat}, which means that there is a model structure on the category of collections transferable to the category of operads (see \cite{SymBat} for the details).

For $n=1$ the relationships between operadic categories above  is simple. The following categories are isomorphic to the category of nonsymmetric operads
$$O_1\simeq LCO_1 \simeq QO_1$$
and we have a classical adjunction between nonsymmetric operads and symmetric operads. All this is true on the level of homotopy categories.

For $n=2$ we have the following  diagram of categories and right and left adjoint functors:

{\unitlength=1mm

\begin{picture}(40,34)(-29,0)

\put(13,25){\makebox(0,0){\mbox{$O_2$}}}
\put(15,21){\vector(1,-1){10}}
\put(12,11){\vector(0,1){10}}
\put(18,24){\vector(1,0){22}}
\put(40,26){\vector(-1,0){22}}
\put(24,27){\shortstack{\mbox{$\scriptstyle Des_2$}}}
\put(24,22){\shortstack{\mbox{$\scriptstyle Sym_2$}}}
\put(45,25){\makebox(0,0){\mbox{$ SO$}}}
\put(44,21){\vector(0,-1){10}}
\put(46,11){\vector(0,1){10}}
\put(13,7){\makebox(0,0){\mbox{$ LCO_2$}}}
\put(25,8){\vector(-1,0){7}}
\put(18,6){\vector(1,0){7}}
\put(41,8){\vector(-1,0){7}}
\put(34,6){\vector(1,0){7}}
\put(45,7){\makebox(0,0){\mbox{$BO$}}}
\put(30,7){\makebox(0,0){\mbox{$QO_2$}}}
\put(8,15){\shortstack{\mbox{$\scriptstyle I_2$}}}
\put(15,15){\shortstack{\mbox{$\scriptstyle L_2$}}}
\put(40,15){\shortstack{\mbox{$\scriptstyle U_2$}}}
\put(47,15){\shortstack{\mbox{$\scriptstyle F_2$}}}
\put(20,3.5){\shortstack{\mbox{$\scriptstyle K_2$}}}
\put(20,9){\shortstack{\mbox{$\scriptstyle J_2$}}}
\put(36,3.5){\shortstack{\mbox{$\scriptstyle B_2$}}}
\put(36,9){\shortstack{\mbox{$\scriptstyle A_2$}}}
\end{picture}}

In this diagram the functor $Des_2$ is right adjoint to $Sym_2$ (see \cite{SymBat,EHBat} for the construction). The functors $I_2$ and $J_2$ are natural  inclusions. The functor $K_2$ is left adjoint to $J_2$ and $L_2$ is left adjoint to the composite $J_2\cdot I_2 .$ Using the theory of internal operads from \cite{EHBat} one can show that  $L_2$  on the level of collections is given by the left Kan extension along the localisation functor $l_2:Q_2\rightarrow \Pi Q_2:$ 
\begin{equation}\label{kan} L_2(A) = Lan_{l_2}(A) .\end{equation} 
We have also the same formula for $K_2 .$  The functor $A_2$ is a right adjoint and $B_2$ is a left adjoint part of the equivalence constructed in the section \ref{vs}. Finally, $U_2$ is the functor which produces a braided operad from a symmetric operad by pulling back along the functor $p:\Br\rightarrow \Sm$ and $F_2$ is its left adjoint given by quotienting with respect to the action of the pure braid groups .

\begin{theorem}\label{lc=2}\begin{itemize}
   
\item The homotopy category of locally constant $2$-operads and the homotopy category of quasisymmetric $2$-operads and the homotopy category of braided operads are equivalent.
\item The functor of symmetrisation $Sym_2 $ can be factorised as $L_2\cdot B_2 \cdot F_2 .$
\item A base space $X$ is a $2$-fold loop space (up to group completion) if and only if it is an algebra of  a contractible $2$-operad , if and only if it is an algebra of a contractible braided operad (Fiedorowicz's recognition principle \cite{fied}).       
\end{itemize}
\end{theorem}

\Proof Since $\Pi Q_2$ is a groupoid, the localisation functor $l_2$ is locally constant in the sense of \cite[1.14]{cis}.
By Formal Serre spectral sequence \cite[Prop. 1.15]{cis} we get 
that the homotopy left Kan extension along $l_2$ is a left adjoint to the restriction functor between homotopy categories of collections. The functor $l_2$ induces a weak equivalence of the nerves and so, by Quillen theorem B, it is also aspherical in the sense of \cite[1.4]{cis}. So, by \cite[Prop. 1.16]{cis} the homotopy left Kan extension along $l_2$ is an equivalence of homotopy categories of collections.    
 
Taking into account  the formula (\ref{kan}) we see that to prove the equivalence of homotopy categories of operads it is enough to show   that 
for an $n$-operad $A$ ($1\le n\le \infty$) there exists  a cofibrant replacement $B(A)$ such that the underlying $Q_n$-collection of $B(A)$ is cofibrant in the projective model structure.   

Recall \cite{SymBat}, that $\ph^n$ is the categorical symmetric operad representing the $2$-functor of internal pruned  $n$-operads. In particular an $n$-operad $A$ is represented by an operadic functor $\tilde{A}:\ph^n \rightarrow V^{\bullet}.$ 
If we forget about operadic structures then for any $k\ge 0 $ we will have a functor $\tilde{A}_k: \ph^n_k \rightarrow V .$ 
Take the bar-resolution  $B(L,L,C(\tilde{A})),$ where
$(L,\mu,\epsilon)$ is the monad on the functor category $[d(\ph^n),V]$   generated by restriction and left Kan extension along the inclusion of discretisation  $d(\ph^n_k)$ of $\ph^n_k$  to $\ph^n_k$ and $C(A)$ is the termwise cofibrant replacement of the underlying $n$-collection of $A .$  These functors    
for all $k\ge 0$ form an operadic functor ${\bf B}(A):\ph^n\rightarrow V^{\bullet}$ and, hence, determine an $n$-operad $B(A)$ which is a cofibrant replacement for $A$ \cite{EHBat}. 

Since ${\bf B}_k(A)$ is a bar-construction on cofibrant collection it is cofibrant in the projective model category of functors. Recall also that 
there is a symmetric categorical operad $\rh^n$  representing the $2$-functor of internal reduced $n$-operads \cite{SymBat} and  a projection $p:\ph^n\rightarrow \rh^n.$ A typical fiber (in a strict sense)  of this projection over an object $w\in \rh^n$  is a category with a terminal object $s(w).$ The map $s$ assembles to the (nonoperadic) functor $s:\rh^n\rightarrow \ph^n ,$ which is by definition a section of $p $ and it is also a right adjoint to $p .$ The counit of this adjunction is the identity and the unit is the unique map to the terminal object $s(w).$

The simple calculations with this adjunction shows that 
the restriction functor $s^*$  preserves the cofibrant objects for projective model structures and so $s^*({\bf B}_k(A))$  is cofibrant.
     
There is also an inclusion $j:J^{op}_n\rightarrow \rh^n $ \cite{SymBat}. 
It is not hard to see also that the categories $J^{op}_n(k)$ and $\rh^n_k$ are Reedy categories. Recall, that the objects of $\rh^n_k$ are planar trees decorated by pruned $n$-trees (i.e. $n$-ordinals). One can choose the total number of edges of $n$-trees in a decorated planar tree as a degree function and see that each morphism decreases strictly this function.\footnotemark\footnotetext{In fact, $\rh^n_k$ is a poset but we did not provide a proof of this fact in \cite{SymBat}.}.    

It follows from these considerations that the functor $s^*({\bf B}_k(A))$ satisfies the following property characterising cofibrant objects in the projective model categories for functor categories over Reedy categories:
\begin{equation}\label{cofinal}colim( s^*({\bf B}_k(A))(w)) \rightarrow s^*({\bf B}_k(A))(T)\end{equation} 
is a cofibration. Here the colimit is taken over the category of all 
$w \rightarrow T , \ w\ne T$ in $\rh^n_k .$ It was proved in \cite{SymBat} that $J^{op}_n(k)$ is cofinal in $\rh^n_k .$ Exactly the same argument shows that in the colimit (\ref{cofinal}) one can replace $w\in rh^n_k$ by the objects from $J^{op}_n(k).$ 
And, therefore the restriction  $j^*s^*({\bf B}_k(A))$ is cofibrant as well. 

The quotient functor $q:J^{op}_n(k)\rightarrow Q^{op}_n(k)$ induces the restriction functor $q^*$ on functor categories which is fully faithful. It follows from this that $q^*$ reflects cofibrations. We observe that 
$$q^*(u(B(A)) = j^*s^*({\bf B}_k(A))$$
and so $u(B(A))$ is cofibrant. Hence the first statement of the theorem is proved.

 The statement about symmetrisation is obvious since $Des_2 = U_2\cdot A_2\cdot J_2\cdot I_2 .$ Finally, a contractible operad is locally constant so the third statement follows from the first statement, Theorem 8.6 from \cite{SymBat} and the fact that the functors 
$U_2,A_2,J_2,I_2$ preserve endomorphism operads.

\Q

For $3\le n\le \infty$ the corresponding diagram is:  

{\unitlength=1mm

\begin{picture}(40,34)(-29,0)
\put(39.5,15){\shortstack{\mbox{$\scriptstyle A_n$}}}
\put(47,15){\shortstack{\mbox{$\scriptstyle B_n$}}}
\put(28,3.5){\shortstack{\mbox{$\scriptstyle K_n$}}}
\put(28,9){\shortstack{\mbox{$\scriptstyle J_n$}}}
\put(8,15){\shortstack{\mbox{$\scriptstyle I_n$}}}
\put(21,15){\shortstack{\mbox{$\scriptstyle L_n$}}}
\put(13,25){\makebox(0,0){\mbox{$O_n$}}}
\put(17,21){\vector(2,-1){21}}
\put(12,11){\vector(0,1){10}}
\put(18,24){\vector(1,0){22}}
\put(40,26){\vector(-1,0){22}}
\put(24,27){\shortstack{\mbox{$\scriptstyle Des_n$}}}
\put(24,22){\shortstack{\mbox{$\scriptstyle Sym_n$}}}
\put(45,25){\makebox(0,0){\mbox{$ SO$}}}
\put(44,21){\vector(0,-1){10}}
\put(46,11){\vector(0,1){10}}
\put(13,7){\makebox(0,0){\mbox{$ LCO_n$}}}
\put(19,8){\vector(1,0){20}}
\put(39,6){\vector(-1,0){20}}
\put(45,7){\makebox(0,0){\mbox{$QO_n$}}}
\end{picture}}

 \begin{theorem}\begin{itemize}
\item For $3\le n \le \infty$ the category of symmetric operads is equivalent to the category of quasisymmetric $n$-operads; 

\item For $3\le n < \infty$ a base space $X$ is an $n$-fold loop space (up to group completion) if and only if it is an algebra of  a contractible $n$-operad;
\item The homotopy category of locally constant $\infty$-operads, the homotopy category of quasisymmetric $\infty$-operads and the homotopy category of symmetric operads are equivalent.
\item  A base space $X$ is an infinite loop space (up to group completion) if and only if it is an algebra of  a contractible $\infty$-operad if and only if it is an algebra of a contractible symmetric operad (May's recognition principle \cite{May}).
\end{itemize}
\end{theorem}

\Proof The proof is analogous to the proof of  Theorem \ref{lc=2}.

\Q

An interesting question which we do not consider here is the existence of model structures on the various categories of operads. The results of \cite{cis} indicate that this might be possible. But it is a subject for a future paper.  

\   

\noindent {\bf Acknowledgements.}  I would like to thank Denis-Charles Cisinski for his nice answers \cite{cis} to my some time naive questions.    I   wish to express my  gratitude to
C.Berger,  I.Galvez, E.Getzler, V.Gorbunov, A.Davydov, R.Street, A.Tonks,  M.Weber for many useful discussions and to my anonymous referee for  useful comments concerning the presentation of the paper.

I also  gratefully acknowledge  the financial
support of Scott Russel Johnson Memorial Foundation, Max Plank
Institut f\"{u}r Mathematik and Australian Research Council (grant
No.~DP0558372). 

\renewcommand{\refname}{Bibliography.}

\end{document}